\documentclass[12pt,twoside]{article}
\usepackage{times}
\usepackage{amsmath,amssymb}
\usepackage{color}
\usepackage{graphicx}

\allowdisplaybreaks

\def\p{\partial}
\def\ve{\varepsilon}
\def\f{\frac}
\def\na{\nabla}
\def\la{\lambda}

\def\al{\alpha}
\def\t{\tilde}
\def\q{\quad}
\def\vp{\varphi}
\def\O{\Omega}
\def\th{\theta}
\def\g{\gamma}
\def\G{\Gamma}
\def\si{\sigma}
\def\dl{\delta}
\def\p{\partial}
\def\ve{\varepsilon}
\def\f{\frac}

\def\na{\nabla}
\def\la{\lambda}
\def\al{\alpha}

\def\t{\tilde}

\def\o{\omega}
\def\O{\Omega}
\def\vp{\varphi}
\def\th{\theta}

\def\g{\gamma}
\def\G{\Gamma}

\def\si{\sigma}
\def\dl{\delta}

\def\q{\qquad}

\def\ds{\displaystyle}
\def\dP{\dot\Phi}
\def\ss{\ds\sum}
 \allowdisplaybreaks

\begin{document}
 \footskip=0pt
 \footnotesep=2pt
\let\oldsection\section
\renewcommand\section{\setcounter{equation}{0}\oldsection}
\renewcommand\thesection{\arabic{section}}
\renewcommand\theequation{\thesection.\arabic{equation}}
\newtheorem{claim}{\noindent Claim}[section]
\newtheorem{theorem}{\noindent Theorem}[section]
\newtheorem{lemma}{\noindent Lemma}[section]
\newtheorem{proposition}{\noindent Proposition}[section]
\newtheorem{definition}{\noindent Definition}[section]
\newtheorem{remark}{\noindent Remark}[section]
\newtheorem{corollary}{\noindent Corollary}[section]
\newtheorem{example}{\noindent Example}[section]

\title{On global multidimensional supersonic flows with vacuum states at infinity}

\author{Xu Gang$^{1}$, \quad Yin
Huicheng$^{2}$\footnote{Xu Gang was supported
 by the National Natural Science Foundation of China (No.11101190)
and Natural Science Fundamental Research Project of Jiangsu Colleges
(No.10KLB110002); Yin Huicheng was
supported by the NSFC (No.~10931007, No.~11025105) and the Priority
Academic Program Development of Jiangsu Higher Education
Institutions.}\vspace{0.5cm}\\
\small 1. Faculty of Science, Jiangsu University, Zhenjiang, Jiangsu
212013, China.\\
\small 2. Department of Mathematics
and IMS, Nanjing University, Nanjing 210093, China.\\
}

\date{}
\maketitle

\centerline {\bf Abstract} \vskip 0.3 true cm

In this paper, we are
concerned with the global existence and  stability of a smooth supersonic flow
with vacuum state at infinity in a 3-D infinitely long divergent nozzle.
The flow is described by a 3-D steady potential equation, which is multi-dimensional quasilinear
hyperbolic (but degenerate at infinity) with respect to the supersonic direction,
and whose linearized part admits the form
$\p_t^2-\ds\f{1}{(1+t)^{2(\g-1)}}(\p_1^2+\p_2^2)+\ds\f{2(\g-1)}{1+t}\p_t$ for $1<\g<2$. From the physical point of view,
due to the expansive geometric property of the divergent nozzle and the mass conservation of gas, the
moving gas in the nozzle will gradually
become rarefactive and tends to a vacuum state at infinity, which implies that such a smooth
supersonic flow should be globally stable for small perturbations since there are no strong
resulting compressions in the motion of the flow. We will confirm such a global stability phenomena
by rigorous mathematical proofs and further show that there do not exist
vacuum domains in any finite part of the nozzle.

\vskip 0.3 true cm

\vskip 0.3 true cm


\vskip 0.4 true cm
\centerline{\bf $\S 1$. Introduction  and main results}
\vskip 0.3 true cm

In this paper, we are concerned with the global existence and stability of
a smooth supersonic polytropic gas with vacuum state at infinity in a 3-D infinitely long divergent nozzle.
The divergent nozzle is described by the domain $\O=\{x=(x_1, x_2, x_3)\in\Bbb R^3: x_1^2+x_2^2\le\tan^2\vp_0 x_3^2,
x_1^2+x_2^2+x_3^2\ge 1, x_3>0\}$
with $\vp_0\in (0, \ds\f{\pi}{2})$ (see the Figure 1 below), and the potential function $\Phi$
of irrotational polytropic gas satisfies the
following steady potential equation in $\O$:
$$
\ds\sum_{i=1}^3((\p_i\Phi)^2-c^2(\rho))\p_i^2\Phi+2 \ds\sum_{1\le i<j\le
3}\p_i\Phi\p_j\Phi\p_{ij}^2\Phi=0,\eqno{(1.1)}
$$
where $\p_i=\p_{x_i}$ ($1\le i\le 3$), $c(\rho)=\sqrt{P'(\rho)}$ is the local sound speed, $P(\rho)$ is the pressure,
$\rho$ is the density, and state equation is given by $P(\rho)=\rho^{\g}$ with $1<\g<2$ (for the air, $\g\approx 1.4$).
Moreover, the density $\rho=\rho(\na_x\Phi)$ can be determined by the Bernoulli's law:
$$\f{1}{2}|\na_x \Phi|^2+\f{\g}{\g-1}\rho^{\g-1}=C_0\equiv \f12 q_0^2+\f{\g}{\g-1}\rho_0^{\g-1},\eqno{(1.2)}$$
where $\na_x=(\p_1, \p_2, \p_3)$, and $q_0>c(\rho_0)$ (this means that the flow at the entrance is supersonic
along the radial direction).
Without loss of generality and for convenience,
$C_0=1$ will be always assumed in the whole paper.

\vskip 0.5 true cm
\includegraphics[width=12cm,height=6.5cm]{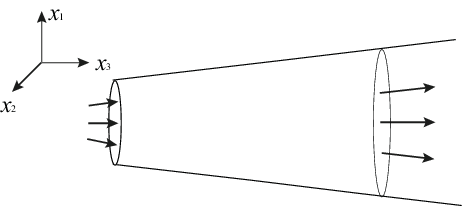}

\centerline{\bf Figure 1. Supersonic flow in a 3-D divergent nozzle}
\vskip 0.8 true cm

Denote the divergent nozzle wall by $\Sigma=\{x:
x_1^2+x_2^2=\tan^2\vp_0 x_3^2, x_1^2+x_2^2+x_3^2\ge 1, x_3>0\}$, then
$\Phi$ satisfies the following fixed boundary condition on $\Sigma$:
$$x_1\p_1\Phi+x_2\p_2\Phi-\tan^2\vp_0x_3\p_3\Phi=0.\eqno{(1.3)}$$

Due to the divergent geometric property of $\O$, it is convenient to work in
the spherical  coordinates $(r, \th, \vp)$:
$$(x_1, x_2, x_3)=(r\cos\th\sin\vp, r\sin\th\sin\vp, r\cos\vp),\eqno{(1.4)}$$
where
$r=\sqrt{x_1^2+x_2^2+x_3^2}$, $0\le\th\le 2\pi$ and $0\le\vp\le\vp_0$.

Under the coordinate transformation (1.4), (1.1) becomes
\begin{align*}
&\big((\p_r\Phi)^2-c^2(\rho)\big)\p_r^2\Phi+
\f{1}{r^2\sin^2\vp}\bigg(\f{1}{r^2\sin^2\vp}(\p_{\th}\Phi)^2-c^2(\rho)\bigg)\p_{\th}^2\Phi
+\f{1}{r^2}\bigg(\f{1}{r^2}(\p_{\vp}\Phi)^2-c^2(\rho)\bigg)\p_{\vp}^2\Phi\\
&\qquad +\f{2\p_r\Phi\p_{\th}\Phi}{r^2\sin^2\vp}\p_{r\th}^2\Phi+\f{2}{r^2}\p_r\Phi\p_{\vp}\Phi\p_{r\vp}^2\Phi
+\f{2\p_{\th}\Phi\p_{\vp}\Phi}{r^4\sin^2\vp}\p_{\th\vp}^2\Phi
-\f{1}{r^3}\bigg(2r^2c^2(\rho)+(\p_{\vp}\Phi)^2\\
&\qquad +\f{1}{\sin^2{\vp}}(\p_{\th}\Phi)^2\bigg)\p_r\Phi
-\f{\cot\vp}{r^4}\bigg(r^2c^2(\rho)+\f{1}{\sin^2\vp}(\p_{\th}\Phi)^2\bigg)\p_{\vp}\Phi=0.\tag{1.5}
\end{align*}

In particular, if the solution $\Phi$ of (1.5) is axially symmetric, namely, $\Phi(r,\th,\vp)\equiv \Phi(r,\vp)$
is independent of the variable $\th$, then (1.5) becomes
\begin{align*}
&\big((\p_r\Phi)^2-c^2(\rho)\big)\p_r^2\Phi
+\f{1}{r^2}\bigg(\f{1}{r^2}(\p_{\vp}\Phi)^2-c^2(\rho)\bigg)\p_{\vp}^2\Phi
+\f{2}{r^2}\p_r\Phi\p_{\vp}\Phi\p_{r\vp}^2\Phi\\
&\qquad -\f{1}{r^3}\big(2r^2c^2(\rho)+(\p_{\vp}\Phi)^2\big)\p_r\Phi
-\f{c^2(\rho)}{r^2}cot\vp\p_{\vp}\Phi=0.\tag{1.6}
\end{align*}

Here we point out that some coefficients in (1.5) or (1.6) admit strong singularities near $\vp=0$.
Consequently, in order to overcome
the difficulties arisen by the singularities
near $\vp=0$, we require to  rewrite (1.5) or (1.6) by introducing some
smooth vector fields tangent to the sphere $\Bbb S^2$ as in [15].

Set
\begin{equation}
\left\{
\begin{aligned}
&Z_1=x_1\p_2-x_2\p_1=\p_{\th},\\
&Z_2=x_2\p_3-x_3\p_2=-\cot\vp\cos\th\p_{\th}-\sin\th\p_{\vp},\\
&Z_3=x_3\p_1-x_1\p_3=-\cot\vp\sin\th\p_{\th}+\cos\th\p_{\vp}.
\end{aligned}
\right.\tag{1.7}
\end{equation}
Then it follows from a direct computation that (1.5) or (1.6) has such a new form
\begin{align*}
&((\p_r\Phi)^2-c^2(\rho))\p_r^2\Phi+\f{2\p_r\Phi}{r^2}\sum\limits_{i=1}^3Z_i\Phi\p_rZ_i\Phi
-\f{c^2(\rho)}{r^2}\sum\limits_{i=1}^3Z_i^2\Phi
+\f{1}{r^4}\sum\limits_{i,j=1}^3Z_i\Phi
Z_j\Phi
Z_iZ_j\Phi\\
&\quad +\sum\limits_{i,j=1}^3\f{C_{ij}(\o)}{r^3}\p_r\Phi Z_i\Phi
Z_j\Phi+\sum\limits_{i,j,k=1}^3\f{C_{ijk}(\o)}{r^4}Z_i\Phi
Z_j\Phi
Z_k\Phi-\f{2c^2(\rho)}{r}\p_r\Phi=0,\tag{1.8}
\end{align*}
where $\o=\ds \f{x}{r}$,
$C_{ij}(\o)=C_{ij}(\ds \f{x}{r})$ and $C_{ijk}(\o)=C_{ijk}(\ds\f{x}{r})$
are smooth functions on their arguments.

Meanwhile, the fixed boundary condition (1.3) can be changed as
$$
x_1Z_3\Phi-x_2Z_2\Phi=0\qquad \text{on}~~\Sigma.\eqno{(1.9)}
$$

Especially, for the axially symmetric solution $\Phi$, the boundary condition on $\Sigma$ is

$$
Z_2\Phi=Z_3\Phi=0\qquad \text{on}~~\Sigma.\eqno{(1.10)}
$$

In addition, we impose the following initial axially symmetric perturbations:
$$\Phi(1,\th, \vp)=\ve\Phi_0(\vp),\qquad
\p_r\Phi(1,\th,\vp)=q_0+\ve\Phi_1(\vp),\eqno{(1.11)}
$$
where $\ve>0$ is a small parameter, and
$\Phi_i(\vp)\in C_0^{\infty}[0,\vp_0)$ ($i=0,1$). In fact, such kinds of initial conditions (1.11)
can be easily realized by small axially symmetric perturbations  on  the initial
density and velocity of irrotational gas.

\vskip 0.5 true cm
\includegraphics[width=11cm,height=6.0cm]{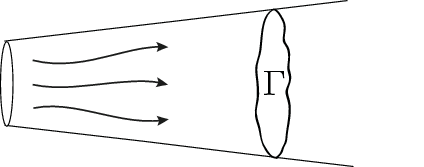}

\centerline{\bf Figure 2. Perturbed supersonic flow in the 3-D divergent nozzle}
\vskip 0.8 true cm

Let $\G=\{r=r(\th, \vp): 0\le\th\le 2\pi, 0\le\vp\le\vp_0\}$ be any $C^1-$smooth
cross section of $\O$ (see the Figure 2 above). Denote the positive constant
$m_{\ve}
=2\pi\int_0^{\vp_0}\rho_0^{\ve}(\vp)
(q_0+\ve\Phi_1(\vp))sin\vp d\vp$, where the initial density $\rho_0^{\ve}(\vp)=(\ds\f{\g-1}{\g})^{\f{1}{\g-1}}\biggl\{\ds\f{\g}{\g-1}\rho_0^{\g-1}-\f12
\biggl(2q_0\ve \Phi_1+\ve^2\Phi_1^2+\ve^2(\Phi'_0)^2\biggr)\biggr\}^{\f{1}{\g-1}}$
is determined by the Bernoulli's law (1.2). The main result in our paper is:

{\bf Theorem 1.1.} {\it There exists a constant $\ve_0>0$ depending on $q_0,\rho_0$
and $\gamma$ such that problem (1.8) with (1.10)-(1.11)
possesses a global $C^{\infty}$ supersonic solution $\Phi(x)$ for
$\ve<\ve_0$ and the mass of gas on any smooth cross surface $\G$ is conserved, namely,
$\int_{\G}\rho\na_x\Phi\cdot {\overrightarrow {n}}dS
\equiv m_{\ve}$, where $\overrightarrow {n}$ stands for the unit outward normal direction of
$\G$. Moreover, $\rho(x)>0$ and $\ds\lim_{r\to\infty}\rho(x)=0$ hold in the whole $\O$.}

{\bf Remark 1.1.} {\it From Theorem 1.1, one easily knows that there do not exist
vacuum domains in any finite part of $\O$ for the problem (1.8) together with (1.10)-(1.11).}

{\bf Remark 1.2.} {\it  For the small arbitrarily (not axially symmetric) perturbed supersonic flow in $\O$,
which is determined by the equation
(1.8) together with (1.9) and the initial data $(\Phi(1,\th,\vp),$ $\p_r\Phi(1,\th,\vp))
=(\ve \Phi_0(\th,\vp), q_0+\ve\Phi_1(\th,\vp))$ with $\Phi_i(\th,\vp)\in C_0^{\infty}([0, 2\pi]\times [0, \vp_0))$ $(i=0,1)$,
we can also solve the global stability problem  as in Theorem 1.1 by analogous but much more complicated analysis.
Nevertheless, due to  the lengthy formulas and too heavy computations, we do not give out the related details of
proof procedure here.}

{\bf Remark 1.3.} {\it  By the same analysis in this paper, Theorem 1.1 can be extended into the curved
2-D or 3-D divergent nozzles with small and arbitrary perturbations of straight boundaries (one can see the
following Figure 3 and Figure 4).}

\vskip 0.3 true cm
\includegraphics[width=11cm,height=5.5cm]{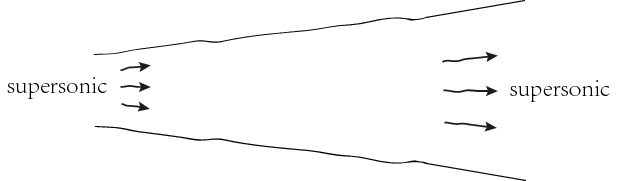}

\centerline{\bf Figure 3. 2-D global smooth supersonic flow in a curved divergent nozzle}

\vskip 0.8 true cm
\includegraphics[width=11cm,height=5.5cm]{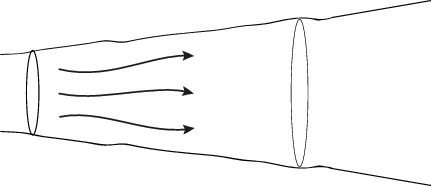}

\centerline{\bf Figure 4. 3-D  global smooth supersonic flow in a curved divergent nozzle}
\vskip 0.8 true cm

{\bf Remark 1.4.} {\it  For the de Laval nozzle, which is constructed  by
a converging ``entry'' section and a diverging ``exhaust''
section, when the supersonic flow is formed across the sonic  curve
in the slowly variable nozzle and the infinite long nozzle walls approach two symmetric lines (see the Figure 5 below),
then our Theorem 1.1 illustrates that the smooth supersonic flow exists globally for the small
perturbed state. On the other hand, if the de Laval nozzle is finitely long and
an appropriately large exit pressure
$p_e$ is given,  as stated in Section 147 of [5],
at a certain place in
the diverging part of the nozzle a shock front intervenes and the gas is compressed and slowed down to
subsonic speed (see the Figure 6 below). This phenomenon has been extensively studied, especially the
stability problem of a transonic shock is completely solved  for
a general class of 2-D de
Laval nozzles whose divergent parts are small and arbitrary perturbations of divergent angular domains for
the full steady compressible Euler system in [16].}

\vskip 0.5 true cm
\includegraphics[width=13cm,height=5.5cm]{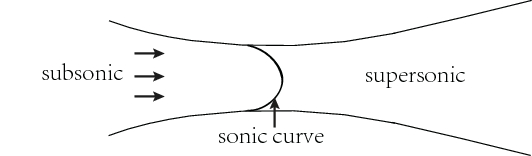}

\centerline{\bf Figure 5. Global continuous transonic flow in an infinite long de Laval nozzle}
\vskip 0.8 true cm

\vskip 0.5 true cm
\includegraphics[width=16cm,height=5.5cm]{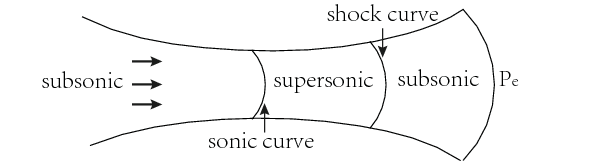}

\centerline{\bf Figure 6. Stability of a transonic shock in a finitely long de Laval nozzle}
\vskip 0.8 true cm

{\bf Remark 1.5.} {\it The nonlinear equation (1.1) in our case is actually a two dimensional
quasilinear degenerate wave equation
if one regards $r$ as the time since the flow is
supersonic in $r$-direction, whose linearized part is like $\p_t^2-\ds\f{1}{(1+t)^{2(\g-1)}}
(\p_1^2+\p_2^2)+\ds\f{2(\g-1)}{1+t}\p_t$ (one can see Remark 3.1 below in $\S 3$). On the other hand,
if we consider the Cauchy initial data problem of (1.1) which is of a small
perturbation with respect to the uniform constant density $\rho_0$ and velocity $(0,0,q_0)$
\begin{equation}
\left\{
\begin{aligned}
&\ds\sum_{i=1}^3((\p_i\Phi)^2-c^2(\rho))\p_i^2\Phi+2 \ds\sum_{1\le i<j\le
3}\p_i\Phi\p_j\Phi\p_{ij}^2\Phi=0, \qquad x_3\ge 0,\\
&\Phi(x)|_{x_3=0}=\ve\Phi_0(x_1,x_2),\quad \p_3\Phi(x)|_{x_3=0}=q_0+\ve\Phi_1(x_1,x_2),\qquad (x_1, x_2)
\in\Bbb R^2,
\end{aligned}
\right.\tag{1.12}
\end{equation}
where $q_0>c(\rho_0)$, and $\Phi_i(x_1,x_2)\in C_0^{\infty}(\Bbb R^2)$ ($i=0,1$), then by a direct verification,
one sees that (1.12) does not fulfill the ``null-condition'' put
forward in [4] and [14]. Therefore, in terms of the extensive
results of [1-2], [9], [21], [28] and so on, the classical solution to
(1.12) will blow up for finite $x_3$. However, compared this blowup result with our Theorem 1.1,
we obtain the global existence of a smooth solution to (1.12) together with the divergent
nozzle wall condition due to the rarefactive property of
supersonic gas.}

{\bf Remark 1.6.} {\it If the initial density contains vacuum, the local well-posedness results of
compressible Euler system have been extensively studied in [3], [6-7], [12-13], [18], [20], [25-27]
and so on. In the general case, such local classical solution will blow up in finite time as shown in [3], [25] and
the references therein. With respect to the problem in our paper, the vacuum only appears at infinity
and the smooth solution exists globally.}

{\bf Remark 1.7.} {\it If the initial velocity $u_0(x)$ of gas forces particles to spread out, roughly speaking,
$u_0(x)$ is close to a linear field, which means $\ds{\overline\lim}_{
|x|\rightarrow \infty}|u_0(x)|
=\infty$, the authors in [11] and [22] have proved the global existence
of smooth solutions to the Cauchy problem of compressible Euler system. Here we emphasize that
our initial data (1.11) are not the cases posed in [11] or [22]
(for example, one can see Theorem 1 of [11]).}

Let us comment on the proof of Theorem 1.1. Since the local solvability
of problem (1.8) together with (1.10)-(1.11) has been known
as long as the vacuum does not appear, we will use the continuous induction method  to
prove Theorem 1.1. To achieve
this objective, we need to establish the global
energy estimates with suitable anisotropic weights for (1.8) with
(1.10)-(1.11), which is degenerate at infinity and admits a linear part as follows:
$\p_r^2-\ds\f{1}{(1+r)^{2(\g-1)}}(\p_1^2+\p_2^2)+\ds\f{2(\g-1)}{1+r}\p_r$. Based on such estimates,
one then obtains the absence of vacuum
for any finite $x_3>0$ in $\O$, the global
existence, stability,  and the asymptotic behavior of the
solution. The key ingredients in
the analysis to obtain weighted energy estimates are to look
for an appropriate multiplier and the suitable anisotropic weights,  derive available boundary conditions of higher order
derivatives of $\Phi$ on the boundary $\Sigma$ and search for the required weighted Sobolev
interpolations. Finding a suitable multiplier and anisotropic weights are not easy due to the
following reasons: Firstly, to obtain the global existence with no vacuum state requires to
establish the estimates independent of $x_3$ and $\na_x^{\al}\Phi$ $(0\le |\al|\le 4)$
on the boundaries as well as in the
interior of the domain $\O$. This leads to strict constraints on
the multiplier and anisotropic weights, as well as makes the computations delicate and
involved. Secondly, as our background solution tends
to vacuum at infinity with different rates for the density and velocity
and their derivatives respectively,
one needs to take some measures to simplify the
coefficients of the nonlinear equation (1.8) so that the procedure to find the
multiplier and anisotropic weights and meanwhile avoid the appearance of vacuum for finite $x_3$ becomes manageable
(one can see more detailed explanations in Remark 2.5 and $\S 5$ below).
Thirdly, the Neumann-type boundary condition (1.10) fulfilled by $\Phi$
arises additional difficulties since there are no enough information
on $\Phi$ itself and its higher order derivatives. Thanks to some delicate analysis
on the radial derivatives and angular derivatives of $\Phi$, which are closely
accompanied by the weighted Sobolev interpolation inequalities in [17], we
finally overcome all these difficulties and obtain a uniform weighted estimate
of $\Phi$ and its higher-order derivatives with no vacuum state for any finite
$x_3>0$ in $\O$.  This eventually establishes Theorem~1.1.

This paper is organized as follows. In \S2, we derive some basic
estimates on the background solution with vacuum at infinity,
and show some preliminary results regarding the weighted Sobolev interpolation
inequalities. In $\S 3$, we reformulate problem (1.8) together with (1.10)-(1.11) by
decomposing its solution as a sum of the background solution and a
small perturbation $\dot\Phi$ so that its linearization can be studied
in a convenient way. In \S4, we will establish a uniform weighted
energy estimate for the corresponding linear problem, where an
appropriate multiplier is constructed. In $\S 5$, the uniform higher-order weighted estimates of
$\dot\Phi$ are established by rather delicate analysis on the radial derivatives and angular derivatives of
$\dot\Phi$, where the domain composition
techniques are applied in order to obtain  the energy estimates of $\dot\Phi$ near $\vp=0$. In \S6, based on
the results in $\S 5$, we complete the proof of Theorem~1.1 by applying Sobolev's embedding theorem and continuous
induction method.

\vskip 0.5 true cm
\centerline{\bf $\S 2$. Background solutions and some preliminaries}
\vskip 0.5 true cm

In this section, at first we analyze the background solution to
(1.6) with (1.10)-(1.11) when the initial data (1.11) are replaced by
$$(\Phi(1, \th, \vp),\quad \p_r\Phi(1, \th, \vp))
=(0, q_0).\eqno{(2.1)}$$
In this case, the density $\rho(x)$ and velocity $u(x)=\na_x\Phi(x)$ in $\O$
have such forms: $\rho(x)=\hat\rho(r)$, $u(x)=\ds\f{x}{r}\hat U(r)$. Consequently, the problem (1.6) with (1.10)
and (2.1) is equivalent to
\begin{equation}
\left\{
\begin{aligned}
&(r^2\hat \rho\hat  U)'(r)=0,\qquad r\ge 1,\\
&\ds\f{1}{2}{\hat U}^2(r)+\f{\gamma}{\gamma-1}{\hat\rho}^{\gamma-1}(r)=1,\qquad r\ge 1,\\
&\hat \rho(1)=\rho_0,\quad \hat U(1)=q_0.
\end{aligned}
\right.\tag{2.2}
\end{equation}

With respect to problem (2.2), we have

{\bf Lemma 2.1.} {\it For $r\ge 1$, (2.2) has a global smooth solution in $\O$ which satisfies
\begin{align*}
&\hat \rho(r)=O(r^{-2})>0,\quad c^2(\hat \rho(r))=O(r^{2(1-\gamma)}),\quad \hat U(r)=\sqrt{2}+O(r^{2(1-\gamma)}),\\
&\hat U'(r)=O(r^{1-2\gamma})>0.\tag{2.3}
\end{align*}
Correspondingly, the potential function $\hat\Phi(r)=\int_1^r\hat U(s)ds$.}

{\bf Remark 2.1.} {\it Lemma 2.1 states an interesting physical phenomenon: along the direction of increasing
area, a supersonic flow is expanded and accelerated, meanwhile becomes more and more rarefactive in the divergent
nozzle. This and more physical phenomena on the supersonic or subsonic
flows in divergent or convergent nozzles can be found in Chapter V of [5].}

{\bf Proof.} It follows from the first equation and the initial data in (2.2) that
$$r^2\hat \rho(r)\hat U(r)=\rho_0q_0.\eqno{(2.4)}$$

This, together with the second equation in (2.2), yields
\begin{equation}
\left\{
\begin{aligned}
&\ds\hat \rho'(r)=-\f{2\rho_0q_0\hat U}{r^3(\hat U^2-c^2(\hat \rho))},\\
&\ds \hat U'(r)=\f{2\hat Uc^2(\hat \rho)}{r(\hat U^2-c^2(\hat \rho))}.
\end{aligned}
\right.\tag{2.5}
\end{equation}
Thus, $\hat \rho'(r)<0$, $\hat U'(r)>0$ and $(\hat U^2-c^2(\hat \rho))'(r)>0$ hold as long as $\hat U^2-c^2(\hat \rho)>0$
and $\hat U>0$. From this, we can also obtain $(\hat U^2-\hat c^2(\hat \rho))(r)\ge (\hat U^2-c^2(\hat \rho))(1)=q_0^2-c^2(\rho_0)>0$
and $\hat U(r)\ge q_0$. On the other hand, if we set $f_1(\hat \rho, \hat U, r)=r^2\hat \rho \hat U-\rho_0q_0$
and $f_2(\hat \rho, \hat U, r)=\ds\f{1}{2}\hat U^2+\f{\gamma}{\gamma-1}\hat \rho^{\gamma-1}-1$, then a direct computation
yields $\ds\f{\p(f_1, f_2)}{\p (\hat \rho, \hat U)}=r^2(\hat U^2-c^2(\hat \rho))\ge q_0^2-c^2(\rho_0)>0$. Thus, $\hat \rho(r)$ and $\hat U(r)$ in (2.2)
exist globally for $r\ge 1$ by implicit function theorem. In addition, (2.3) can be directly obtained by (2.4) and the second
equality in (2.2), and (2.5) respectively. \qquad\qquad \qquad \qquad \qquad \qquad \qquad \qquad
\qquad \qquad \qquad \qquad  \quad $\square$

Next, we cite an important weighted Sobolev interpolation inequality in [17], which will be applied to prove some
crucial weighted inequalities listed in Lemma 2.6 below.

{\bf Lemma 2.2. (see [17])} {Suppose $s, \tau, p, \al, \beta, q, a$ are real
numbers, and $j\ge0, m>0$ are integers, satisfying
\begin{equation}
\left\{
\begin{aligned}
&\ds p,q\ge1, \f{j}{m}\le a\le1, s>0,\\
&\ds \f{1}{s}+\f{\tau}{n}>0,
\f{1}{p}+\f{\al}{n}>0, \f{1}{q}+\f{\beta}{n}>0,\\
&\ds m-j-\f{n}{p} \quad \text {is not a nonnegative integer.}
\end{aligned}
\right.\tag{2.6}
\end{equation}

There exists a positive constant $C$ such that the following
inequality holds for all $v\in C_0^{\infty}(\Bbb R^n)$:
$$
\big||x|^{\tau}\na_x^jv\big|_{L^s}\le
C\big||x|^{\al}\na_x^mv\big|_{L^p}^a\big||x|^{\beta}v\big|_{L^q}^{1-a},\eqno{(2.7)}
$$

if and only if the following conditions hold:
$$
\f{1}{s}+\f{\tau-j}{n}=a(\f{1}{p}+\f{\al-m}{n})+(1-a)(\f{1}{q}+\f{\beta}{n})\quad\text{with
$\tau\le a\al+(1-a)\beta;$}\eqno{(2.8)}
$$

if $\ds\f{1}{q}+\f{\beta}{n}=\f{1}{p}+\f{\al-m}{n}$, then
$$
a(\al-m)+(1-a)\beta+j\le\tau; \eqno{(2.9)}
$$

if $\ds a=\f{j}{m}$, then
$$\tau=a\al+(1-a)\beta.\eqno{(2.10)}$$}

{\bf Corollary 2.3.} {\it For the domain $\O$ defined in $\S 1.1$,  if
$u\in C^m(\bar\Omega)$ and
$$u|_{r\ge T}\equiv 0,\eqno{(2.11)}$$
where $T>1$ is a constant, then we have

(i) (2.7) still holds under the restrictions (2.6)(2.8)-(2.10),
moreover the constant $C$ in the right hand side of (2.7) does not
depend on $T$.

(ii) for $m=2$, $1<\g<2$, $\sigma>0$ and $0<\dl<4\g$,
\begin{align*}
&|r^{\f{2\g+2\sigma-1}{4}}\na_x u|_{L^4(\Omega)}\le C|r^{\f{2\g-1}{2}}\na_x^2u|_{L^2(\Omega)}^{\f{1}{2}}|r^{\sigma}u|_{L^{\infty}(\Omega)}^{\f{1}{2}},\tag{2.12}\\
&|r^{\f{4\g+1}{4}}\na_x u|_{L^4(\Omega)}\le C|r^{\f{2\g+1}{2}}\na_x^2u|_{L^2(\Omega)}^{\f{1}{2}}|r^{\g}u|_{L^{\infty}(\Omega)}^{\f{1}{2}},\tag{2.13}\\
&|r^{\f{8\g-7-\dl}{4}}\na_x u|_{L^4(\Omega)}\le
C|r^{\f{4\g-3-\dl}{2}}\na_x^2u|_{L^2(\Omega)}^{\f{1}{2}}|r^{2(\g-1)}u|_{L^{\infty}(\Omega)}^{\f{1}{2}},\tag{2.14}\\
&|r^{\f{8\g-3-\dl}{4}}\na_x u|_{L^4(\Omega)}\le C|r^{\f{4\g-1-\dl}{2}}\na_x^2u|_{L^2(\Omega)}^{\f{1}{2}}
|r^{2\g-1}u|_{L^{\infty}(\Omega)}^{\f{1}{2}},\tag{2.15}
\end{align*}
where the generic positive constant $C$ is independent of $T$.}

{\bf Proof.} (i) The proof is completely parallel to that of Lemma 2.2
(one can check the details in [17]),
then we omit it here.

(ii) In (2.6)(2.8)-(2.10) of Lemma 2.2, set $s=4, p=2, q=\infty,
a=\f12$ and $j=1, m=2$, one then concludes that:

(2.12) and (2.13) come from (2.7) and the choices of $\tau=\ds\f{2\g+2\si-1}{4}, \al=\ds\f{2\g-1}{2}, \beta=\si$
and $\tau=\ds\f{4\g+1}{4}, \al=\ds\f{2\g+1}{2}, \beta=\g$ respectively;

(2.14) and (2.15) are derived from (2.7) by choosing  $\tau=\ds\f{8\g-7-\dl}{4}, \al=\ds\f{4\g-3-\dl}{2}, \beta=2(\g-1)$
and $\tau=\ds\f{8\g-3-\dl}{4}, \al=\ds\f{4\g-1-\dl}{2}, \beta=2\g-1$ respectively.
\qquad\qquad\qquad\quad $\square$

In order to apply Lemma 2.2 or Corollary 2.3 to derive some weighted Sobolev inequalities in $\O$ without
the restriction (2.11), we require
to establish an extension result as follows:

{\bf Lemma 2.4.} {\it Set $D_T=\{(r,\th,\vp): 1<r<T, 0\le\th\le 2\pi, 0\le\vp<\vp_0\}$
for $T>1$. If $u(x)\in C^{3}(\bar D_T)$ and $r^{\beta}\p^{\alpha}_xu\in L^2(D_T)$ $(|\al|\le 3)$ with some $\beta\in\Bbb R$,
then there exists an extension $Eu\in C^3(\bar D_{\f{9}{8}T})$  of $u$ such that $Eu=u$ in $D_T$,
$Eu|_{r\ge \f{9}{8}T}\equiv 0$ and
$$
|r^{\beta}Eu|_{L^\infty(D_{\f{9}{8}T})}\le C|r^{\beta}u|_{L^\infty(D_{T})},\quad |r^{\beta}\na^{\al}_xEu|_{L^2(D_{\f{9}{8}T})}\le C \ds\sum_{|\nu|\le |\al|}|r^{\beta-|\al|+|\nu|}\na^{\nu}_xu|_{L^2(D_T)},\eqno{(2.16)}
$$
where $C>0$ is independent of $T$.}

{\bf Proof.} In terms of the geometric property of $D_T$, it is convenient to use the spherical coordinate
to work. Denote by $\t u(r,\th,\vp)=u(r cos\th sin\vp , r sin\th sin\vp, r cos\vp)$.
Let $\t E$ be an extension operator defined as follows:
\begin{equation*}
(\t Eu)(r,\th,\vp)=\left\{
\begin{aligned}
&\t u(r,\th,\vp),\q\q\q\q\q 1\le r\le T,\\
&\ds\sum_{j=1}^4\la_j\t u(T+j(T-r),\th,\vp),\q T<r\le \f{9}{8}T\\
\end{aligned}
\right.
\end{equation*}
where
$\ds\sum_{j=1}^4(-j)^k\la_j=1$ for $k=0, 1, 2, 3$.

Noticing that
$$1\le \f{r}{T+j(T-r)}\le\f{9}{4}\qquad \text{for $T\le r\le \f{9}{8}T$ and $0\le j\le 3$},$$
then a direct computation yields
$$|r^{\beta}\t Eu|_{L^\infty(D_{\f{9}{8}T})}\le C|r^{\beta}u|_{L^\infty(D_{T})}
$$
and
\begin{align*}
&|r^{\beta}\nabla_x \t Eu|_{L^2(D_{\f{9}{8}T})}^2\le |r^{\beta}\nabla_x u|_{L^2(D_T)}^2
+|r^{\beta}\na_x \t Eu|_{L^2(D_{\f{9}{8}T}\setminus D_T)}^2\\
&\le |r^{\beta}\nabla_x u|_{L^2(D_T)}^2+|r^{\beta}((\p_r\t Eu)^2+\f{1}{r^2\sin^2\vp}(\p_{\th}\t Eu)^2+\f{1}{r^2}(\p_{\vp}\t Eu)^2)|_{L^2(D_{\f{9}{8}T}\setminus D_T)}^2\\
&\le |r^{\beta}\nabla_x u|_{L^2(D_T)}^2+C|r^{\beta}((\p_ru)^2+\f{1}{r^2\sin^2\vp}(\p_{\th}u)^2+\f{1}{r^2}(\p_{\vp}u)^2)|_{L^2(D_{T}\setminus D_{\f{1}{2}T})}^2\\
&\le C|r^{\beta}\nabla_x u|_{L^2(D_T)}^2.
\end{align*}

Analogously, we have for $|\al|\le 3$
$$|r^{\beta}\na_x^{\al}\t Eu|_{L^2(D_{\f{9}{8}T})}\le C |r^{\beta}\na_x^{\al}u|_{L^2(D_T)}.$$

Choosing a $C^{\infty}-$smooth function $\eta(s)$ with $\eta(s)\equiv 1$
for $s\le 1$ and $\eta(s)\equiv 0$ for $s\ge\ds\f98$ and setting
$$Eu(x)=\eta(\f{r}{T})\t Eu,$$
then $Eu$ satisfies (2.16) and the proof of Lemma 2.4 is completed.\qquad\qquad
\qquad \qquad \quad  $\square$

{\bf Remark 2.2.} {\it From Lemma 2.4, we easily know that Corollary 2.3 still holds when the assumption (2.11)
is removed.}

With respect to the $Z-$fileds introduced in (1.7), we have the following properties by direct verifications
as in [15].

{\bf Lemma 2.5.} {\it
\begin{align*}
&(i)\quad  [Z_1, Z_2]=Z_3, [Z_2, Z_3]=Z_1, [Z_3,Z_1]=Z_2.\\
&(ii)\quad [Z_i, \p_r]=0,  Z_ir=0.\\
&(iii)\quad \ds\na_x f\cdot\na_x g=\p_r f\cdot\p_r
g+\f{1}{r^2}\ds\sum_{i=1}^{3}Z_i f\cdot Z_i g\quad \text{for any $C^1$ smooth
functions $f$ and $g$}.\\
&(iv)\quad  |Z v|\leq r|\na_x v|\quad\text{for any $C^1$ smooth function $v$},
\text{here and below $Z\in\{Z_1, Z_2, Z_3\}$}.\\
&(v)\quad \ds\p_1=\f{x_1}{r}\p_r+\f{x_2}{r^2}Z_1-\f{x_3}{r^2}Z_3;
\quad \p_2=\f{x_2}{r}\p_r+\f{x_3}{r^2}Z_2-\f{x_1}{r^2}Z_1;\quad
\p_3=\f{x_3}{r}\p_r+\f{x_1}{r^2}Z_3-\f{x_2}{r^2}Z_2.
\end{align*}}

{\bf Remark 2.3.} {\it If $u\in C^m(\Bbb R^3)$ with $m\in\Bbb N$, then by Lemma 2.5 we have
$|\na_x^mu|\sim |\p_r^mu|+\ds\f{|\p_r^{m-1}Zu|}{r}+\ds\f{|\p_r^{m-2}Z^2u|}{r^2}+...+\ds\f{|Z^mu|}{r^m}$.}

As direct applications of Remark 2.2 and Lemma 2.5, we have the following inequalities
which will be used again and again in $\S 5$ below.

{\bf Lemma 2.6.} {If $1<\g<2$, $\si\ge\g-1$, $0<\dl<4\g$, $u(x)\in C^4(\bar D_T)$,
then there exists a generic positive constant $C$ independent of $T$ such that
\begin{align*}
(i)\quad &|r^{\f{2\g+2\si-9}{4}}Z^2u|_{L^4(D_T)}\le
C\bigg(\ds\sum_{k=0}^2|r^{\f{2\g-1}{2}-k}\na_x^{2-k}(\f{1}{r}Zu)|_{L^2(D_T)}^{\f{1}{2}}\bigg)
|r^{\si-1}Zu|_{L^{\infty}(D_T)}^{\f{1}{2}}.
\tag{2.17}\\
(ii)\quad  &|r^{\f{4\g-11}{4}}Z^3u|_{L^4(D_T)}\le
C\bigg(\ds\sum_{k=0}^2|r^{\f{2\g+1}{2}-k}\na_x^{3-k}(\f{1}{r}Zu)|_{L^2(D_T)}^{\f{1}{2}}\bigg)
|r^{\g}\na_x (\f{1}{r}Zu)|_{L^{\infty}(D_T)}^{\f{1}{2}}\\
&\q +C\bigg(\ds\sum_{k=0}^2|r^{\f{2\g-1}{2}-k}\na_x^{2-k}(\f{1}{r}Zu)|_{L^2(D_T)}^{\f{1}{2}}\bigg)
|r^{\si-1}Zu|_{L^{\infty}(D_T)}^{\f{1}{2}}.\tag{2.18}\\
(iii)\quad  &|r^{\f{8\g-11-\dl}{4}}\p_rZ^2u|_{L^4(D_T)}+|r^{\f{8\g-7-\dl}{4}}\p_r^2Zu|_{L^4(D_T)}\\
&\quad
\le C\bigg(\ds\sum_{k=0}^2|r^{\f{4\g-3-\dl}{2}-k}\na_x^{2-k}\p_rZu|_{L^2(D_T)}^{\f{1}{2}}\bigg)
|r^{2(\g-1)}\p_rZu|_{L^{\infty}(D_T)}^{\f{1}{2}}.\tag{2.19}\\
(iv)\quad  &|r^{\f{8\g-3-\dl}{4}}\p_r^3u|_{L^4(D_T)}\le
C\bigg(\ds\sum_{k=0}^2|r^{\f{4\g-1-\dl}{2}-k}\na_x^{2-k}\p_r^2u|_{L^2(D_T)}^{\f{1}{2}}\bigg)
|r^{2\g-1}\p_r^2u|_{L^{\infty}(D_T)}^{\f{1}{2}}.\tag{2.20}
\end{align*}}

{\bf Proof.} Let $E$ be the extension operator given in Lemma 2.4, then we have

(i)
\begin{align*}
&|r^{\f{2\g+2\si-9}{4}}Z^2u|_{L^4(D_T)}\\
\le & C|r^{\f{2\g+2\si-1}{4}}\na_x(\f{1}{r}Zu)|_{L^4(D_T)} \qquad \text{(Applying $\f{Z^2}{r^2}=\f{Z}{r}(\f{Z}{r})$
due to Lemma 2.5 (ii))}\\
\le & C|r^{\f{2\g+2\si-1}{4}}\na_xE(\f{1}{r}Zu)|_{L^4(\O)}\\
\le & C|r^{\f{2\g-1}{2}}\na_x^2E(\f{1}{r}Zu)|_{L^2(\O)}^{\f{1}{2}}|r^{\si}E(\f{1}{r}Zu)|_{L^{\infty}(\O)}^{\f{1}{2}}
\quad \text{(Applying (2.12) for $E(\f{1}{r}Zu)$)}\\
\le &
C\bigg(\ds\sum_{k=0}^2|r^{\f{2\g-1}{2}-k}\na_x^{2-k}(\f{1}{r}Zu)|_{L^2(D_T)}^{\f{1}{2}}\bigg)
|r^{\si-1}Zu|_{L^{\infty}(D_T)}^{\f{1}{2}}. \quad\text{(By Lemma 2.4)}\\
\end{align*}

(ii)
\begin{align*}
&|r^{\f{4\g-11}{4}}Z^3u|_{L^4(D_T)}\\
\le &C|r^{\f{4\g+1}{4}}\na_x^2(\f{1}{r}Zu)|_{L^4(D_T)}+C|r^{\f{4\g-3}{4}}\na_x(\f{1}{r}Zu)|_{L^4(D_T)}\\
\le &C|r^{\f{4\g+1}{4}}\na_xE(\na_x(\f{1}{r}Zu))|_{L^4(\O)}+ C|r^{\f{2\g+2\si-1}{4}}\na_x(\f{1}{r}Zu)|_{L^4(D_T)}
\quad \text{(By $\si\ge \g-1$)}\\
\le &C|r^{\f{2\g+1}{2}}\na_x^2E(\na_x(\f{1}{r}Zu))|_{L^2(\O)}^{\f{1}{2}}|r^{\g}E(\na_x (\f{1}{r}Zu))|_{L^{\infty}(\O)}^{\f{1}{2}}\\
&\qquad + C|r^{\f{2\g+2\si-1}{4}}\na_x(\f{1}{r}Zu)|_{L^4(D_T)}\qquad\qquad\qquad  \text{(Applying (2.13) for $E(\na_x(\f{1}{r}Zu))$)}\\
\le &C\bigg(\ds\sum_{k=0}^2|r^{\f{2\g+1}{2}-k}\na_x^{3-k}(\f{1}{r}Zu)|_{L^2(D_T)}^{\f{1}{2}}\bigg)
|r^{\g}\na_x(\f{1}{r}Zu)|_{L^{\infty}(D_T)}^{\f{1}{2}}\\
&+C\bigg(\ds\sum_{k=0}^2|r^{\f{2\g-1}{2}-k}\na_x^{2-k}(\f{1}{r}Zu)|_{L^2(D_T)}^{\f{1}{2}}\bigg)
|r^{\si-1}Zu|_{L^{\infty}(D_T)}^{\f{1}{2}}.\quad\text{(By Lemma 2.4 and (i))}\\
\end{align*}

(iii)
\begin{align*}
&|r^{\f{8\g-11-\dl}{4}}\p_rZ^2u|_{L^4(D_T)}
+|r^{\f{8\g-7-\dl}{4}}\p_r^2Zu|_{L^4(D_T)}\\
\le &C |r^{\f{8\g-7-\dl}{4}}\na_x(\p_rZu)|_{L^4(D_T)}\\
\le &C |r^{\f{8\g-7-\dl}{4}}\na_xE(\p_rZu)|_{L^4(\O)}\\
\le & C
|r^{\f{4\g-3-\dl}{2}}\na_x^2E(\p_rZu)|_{L^2(\O)}^{\f{1}{2}}|r^{2(\g-1)}E(\p_rZu)|_{L^{\infty}(\O)}^{\f{1}{2}}
\quad \text{(Applying (2.14) for $E(\p_rZu)$)}\\
\le & C\bigg(\ds\sum_{k=0}^2|r^{\f{4\g-3-\dl}{2}-k}\na_x^{2-k}\p_rZu|_{L^2(D_T)}^{\f{1}{2}}\bigg)
|r^{2(\g-1)}\p_rZu|_{L^{\infty}(D_T)}^{\f{1}{2}}.\quad\text{(By Lemma 2.4)}\\
\end{align*}

(iv)
\begin{align*}
&|r^{\f{8\g-3-\dl}{4}}\p_r^3u|_{L^4(D_T)}\\
\le & C|r^{\f{8\g-3-\dl}{4}}\na_x(\p_r^2u)|_{L^4(D_T)}\\
\le & C|r^{\f{8\g-3-\dl}{4}}\na_xE(\p_r^2u)|_{L^4(\O)}\\
\le &
C|r^{\f{4\g-1-\dl}{2}}\na_x^2E(\p_r^2u)|_{L^2(\O)}^{\f{1}{2}}|r^{2\g-1}E(\p_r^2u)|_{L^{\infty}(\O)}^{\f{1}{2}}\quad \text{(Applying (2.15) for $E(\p_r^2u)$)}\\
\le &
C\bigg(\ds\sum_{k=0}^2|r^{\f{4\g-1-\dl}{2}-k}\na_x^{2-k}\p_r^2u|_{L^2(D_T)}^{\f{1}{2}}\bigg)
|r^{2\g-1}\p_r^2u|_{L^{\infty}(D_T)}^{\f{1}{2}}.\quad\text{(By Lemma 2.4)}
\end{align*}

Therefore, we complete the proof of Lemma 2.6.\qquad\qquad\qquad\qquad\qquad\qquad\qquad $\square$

Based on Lemma 2.6, we further have

{\bf Lemma 2.7.} {\it If $1<\g<2$, $\si=\min\{1, 2(\g-1)\}$, $0<\dl<4\g$, $u(x)\in C^4(\bar D_T)$, and the following assumptions hold for some constant $M>0$
$$\sum\limits_{0\le
l_1+l_2\le1} r^{l_1}|\p_r^{l_1}Z^{l_2}\p_ru|\le M\ve
r^{-2(\g-1)},\q r^{-1}|Zu|\le M\ve r^{-\si},\q
r^{-1}|Z^2u|\le
M\ve r^{-(\g-1)},\eqno{(2.21)}$$
then
$$
|r^{\f{2\g+2\si-9}{4}}Z^2u|_{L^4}\le  C(M)\ve^{\f{1}{2}}\bigg(\ss_{l=0}^2\big(|r^{\f{4\g-7-\dl+2l}{2}}\nabla_x^l\p_ru|_{L^2}
+|r^{\f{2\g-5+2l}{2}}\nabla_x^l(\f{1}{r}Zu)|_{L^2}\big)\bigg)^{\f{1}{2}},\eqno{(2.22)}$$
and
\begin{align*}
&|r^{\f{4\g-11}{4}}Z^3u|_{L^4}
+|r^{\f{8\g-11-\dl}{4}}\p_rZ^2u|_{L^4}+|r^{\f{8\g-7-\dl}{4}}\p_r^2Zu|_{L^4}
+|r^{\f{8\g-3-\dl}{4}}\p_r^3u|_{L^4}\\
&\le C(M)\ve^{\f{1}{2}}\bigg(\ss_{l=0}^3\big(|r^{\f{4\g-7-\dl+2l}{2}}\nabla_x^l\p_ru|_{L^2}
+|r^{\f{2\g-5+2l}{2}}\nabla_x^l(\f{1}{r}Zu)|_{L^2}\big)\bigg)^{\f{1}{2}},\tag{2.23}
\end{align*}
where $C(M)>0$ is a constant depending on $M$.}

{\bf Remark 2.4.} {\it By $1<\g<2$ and $\si=\min\{1, 2(\g-1)\}$, we can easily conclude
$\si\ge\g-1$, which means that $\si$  satisfies the requirement in Lemma 2.6.}

{\bf Remark 2.5.} {\it (2.21) actually comes from the induction assumptions in Theorem 5.1 on $\dot\Phi$,
where $\dot\Phi$ is the difference between the solution $\Phi$ of (1.8) and the background solution $\hat\Phi$. By (2.21),
we know that $|\p_ru|\le M\ve r^{-2(\g-1)}$ but $|\f{Zu}{r}|\le M\ve r^{-\si}$, and
$|\p_rZu|\le M\ve r^{-2(\g-1)}$ but $|\f{Z^2u}{r}|\le M\ve r^{-(\g-1)}$, which implies that the decay rates
of the radial derivatives and angular derivatives of $u$ are different. Consequently, in order to obtain the
anisotropic energy estimates of $\dot\Phi$ in $\S 5$, we have to pay much attentions on
distinguishing  the different roles
of $\p_r\dot\Phi$ and $Z\dot\Phi$, and this leads to rather involved and delicate analysis.}

{\bf Proof.} In order to prove (2.22)-(2.23), we only verify
$|r^{\f{2\g+2\si-9}{4}}Z^2u|_{L^4}$ to satisfy (2.22)
since the terms in the left hand side of (2.23) can be analogously
done.

It follows from Lemma 2.5, the assumptions on $\g$ and $\dl$, and a direct computation that
\begin{align*}
&\ds\sum_{k=0}^2|r^{\f{2\g-1}{2}-k}\na_x^{2-k}(\f{1}{r}Zu)|_{L^2(D_T)}\\
&\le C\ss_{l=0}^2\bigg(|r^{\f{4\g-7-\dl+2l}{2}}(\nabla_x^l\p_ru)|_{L^2(D_T)}
+|r^{\f{2\g-5+2l}{2}}(\nabla_x^l(\f{1}{r}Zu))|_{L^2(D_T)}\bigg).
\end{align*}

On the other hand, by (2.21) we have
$$|r^{\si-1}Zu|_{L^{\infty}(D_T)}\le
M\ve.
$$

Consequently, by Lemma 2.6 (i), we know that (2.22) holds for $|r^{\f{2\g+2\si-9}{4}}Z^2u|_{L^4}$,
and then the proof of (2.23) can be completed similarly.\qquad\qquad \qquad \qquad \qquad \qquad \qquad
\qquad $\square$

\vskip 0.5 true cm
\centerline{\bf $\S 3$. Reformulation of the problem (1.8) with (1.10)-(1.11)}
\vskip 0.5 true cm

At first, we state a local solvability result on the problem (1.8) with (1.10)-(1.11).

{\bf Lemma 3.1.} {\it There exists a $T_0>1$ such that the problem
(1.8) with (1.10)-(1.11) possesses a local $C^{\infty}$ solution
$\Phi(r,\vp)$ in $\Omega_{T_0}=\{(r,\th,\vp): 1\le r\le
T_0, 0\le\th\le 2\pi, 0\le\vp\le\vp_0\}$. Moreover,
for any $k\in\Bbb N\cup\{0\}$, there exists a positive constant $C_k$
such that
$$||\Phi(r,\th,\vp)-\hat{\Phi}(r)||_{C^k(\Omega_{T_0})}\le C_k\ve,$$
where $\hat\Phi(r)$ is given in Lemma 2.1.}

{\bf Proof.} The quasilinear equation (1.8) is
strictly hyperbolic with respect to the $r-$direction by $\p_r\Phi>c(\rho)$. Thus, by the standard
Picard iteration  as in [19], one can derive that Lemma 3.1 holds.
\qquad \qquad\quad\quad $\square$

Next, we reformulate (1.8) with (1.10)-(1.11).

Let
$\dot{\Phi}=\Phi-\hat{\Phi}$. Then it follows from a direct computation that (1.8) can be reduced to:
$$
\mathcal {L}\dot{\Phi}=\dot
f\quad\text{in $\O$},\eqno{(3.1)}$$
where
\begin{equation}
\left\{
\begin{aligned}
&\mathcal {L}\dot\Phi=\p_r^2\dot\Phi-\ds\f{P_1(r)}{r^2}\sum\limits_{i=1}^3Z_i^2\dot\Phi+
\ds\f{P_2(r)}{r}\p_r\dot\Phi,\\
&\dot f=f_{00}\p_r^2\dot\Phi+\ds\f{1}{r^2}\sum\limits_{i,j=1}^3f_{ij}Z_iZ_j\dot\Phi
+\ds\f{1}{r}\sum\limits_{i=1}^3f_{0i}\p_rZ_i\dot\Phi+f_0\\
\end{aligned}
\right.
\end{equation}
with
\begin{equation}
\left\{
\begin{aligned}
&P_1(r)=\f{\ds{c}^2(\hat\rho)}{\ds\hat {U}^2-{c}^2(\hat\rho)},\\
&P_2(r)=\f{\ds 2}{\ds( \hat {U}^2-{c}^2(\hat\rho))^2}
\big((\g-1)\hat{U}^4+{c}^4(\hat\rho)+\hat{U}^2{c}^2(\hat\rho)\big)\\
\end{aligned}
\right.\tag{3.2}
\end{equation}
and
\begin{equation}
\left\{
\begin{aligned}
&f_{00}=-\f{\ds 1}{\ds \hat{U}^2-
c^2(\hat\rho)}\biggl((\g+1)\hat{U}\p_r\dot\Phi+\ds\f{\g+1}{2}(\p_r\dot\Phi)^2
+\ds\f{\g-1}{2r^2}\sum\limits_{i=1}^3(Z_i\dot\Phi)^2\biggr),\\
&f_{ii}=\ds\f{1}{\hat{U}^2-{c}^2(\hat\rho)}\biggl(\f{\g-1}{2}(\p_r\dot\Phi)^2-(\g-1)\hat{U}\p_r\dot\Phi
+\ds\f{\g-1}{2r^2}\sum\limits_{k=1}^3(Z_k\dot\Phi)^2-\ds\f{1}{r^2}(Z_i\dot\Phi)^2\biggr),\quad
1\le i\le 3,\\
&f_{ij}=f_{ji}=-\ds\f{1}{r^2(\hat{U}^2-{c}^2(\hat\rho))}Z_i\dot\Phi
Z_j\dot\Phi,\q 1\le i\neq j\le 3,\\
&f_{0i}=-\ds\f{1}{r(\hat{U}^2-{c}^2(\hat\rho))}(\hat{U}+\p_r\dot\Phi)Z_i\dot\Phi,\q
1\le i\le 3,\\
&f_0=f_0^1+f_0^2,\\
&\quad f_0^1=\ds\f{1}{\hat{U}^2-{c}^2(\hat\rho)}\biggl\{-\sum\limits_{i,j=1}^3\ds\f{C_{ij}}{r^3}\p_r\dot\Phi
Z_i\dot\Phi
Z_j\dot\Phi-\sum\limits_{i,j,k=1}^3\ds\f{C_{ijk}}{r^4}Z_i\dot\Phi
Z_j\dot\Phi Z_k\dot\Phi\\
&\qquad\qquad
+\ds\f{2}{r}\bigg(\f{\g-1}{2}(\p_r\dot\Phi)^3-(\g-1)\hat{U}(\p_r\dot\Phi)^2
+\ds\f{\g-1}{2r^2}\sum\limits_{i=1}^3(Z_i\dot\Phi)^2\p_r\dot\Phi
+\ds\f{\g-1}{2}\hat{U}(\p_r\dot\Phi)^2\bigg)\biggr\},\\
&\quad f_0^2=\ds\f{\hat U}{\hat{U}^2-{c}^2(\hat\rho)}\biggl(-\sum\limits_{i,j=1}^3\ds\f{C_{ij}}{r^3}Z_i\dot\Phi
Z_j\dot\Phi+\ds\f{\g-1}{r^3}\sum\limits_{i=1}^3(Z_i\dot\Phi)^2\biggr).\\
\end{aligned}
\right.\tag{3.3}
\end{equation}
Here we point out that the terms $f_0^1=
O(\ds\f{(\p_r\dot\Phi Z\dot\Phi)^2}{r^3})+O(\ds\f{(Z\dot\Phi)^3}{r^4})+O(\f{(\p_r\dot\Phi)^2}{r})$ and $f_0^2=O(\ds\f{(Z\dot\Phi)^2}{r^3})$
appeared in $f_0$ will be treated differently
since only such kinds of estimates of $|\p_r\dot\Phi|\le C\ve r^{-2(\g-1)}\to 0$ and $|Z\dot\Phi|\le C\ve r^{1-\si}\not\to
0$ as $r\to\infty$ are derived in $\S 5$ (one can the details in Lemma 5.4 and Lemma 5.5 below). In fact,
$f_0^1$ can be easily estimated since it admits better decay rate with respect to large $r$.

On the nozzle wall $\vp=\vp_0$, $\dot\Phi$ satisfies
$$Z_2\dot\Phi=Z_3\dot\Phi=0.\eqno{(3.4)}$$

In addition, we have the following initial data of $\dP$ from (1.11)
$$\dot\Phi(1,\vp)=\ve\Phi_0(\vp),\qquad
\p_r\dot\Phi(1,\vp)=\ve\Phi_1(\vp),\eqno{(3.5)}
$$

By using Lemma 2.1 and direct computations, we can obtain the following estimates on the coefficients
of $\mathcal {L}\dot\Phi$ in (3.1):

{\bf Lemma 3.1.}
\begin{align*}
&P_1(r)=O(r^{2(1-\g)})>0,\\
&P_2(r)=2(\g-1)+O(r^{2(1-\g)})>0,\\
&P'_1(r)=-\ds\f{c^2(\hat\rho(r))}{r(\hat U^2-c^2(\hat\rho(r))}\bigl(2(\g-1)\hat U^4+2\hat U^2c^2(\hat\rho(r))\bigr)=O(r^{1-2\g})<0,\\
&\ds\f{P'_1(r)}{P_1}=O(r^{-1})<0,\\
&P'_2(r)=\f{c^2(\hat{\rho}(r))}{r(\hat{U}^2-c^2(\hat{\rho}(r)))^3}\bigg(\bigl(12(\g-1)\hat{U}^4
-8(\g-2)\hat{U}^2c^2(\hat{\rho}(r))\bigr)\bigl(\hat{U}^2-c^2(\hat{\rho}(r))\bigr)\\ &\qquad\quad-2(\g+1)\hat{U}^2\bigl(2(\g-1)\hat{U}^4+2\hat{U}^2c^2(\hat{\rho}(r))+2c^4(\hat{\rho}(r))\bigr)\bigg)=O(r^{1-2\g}).
\end{align*}

{\bf Remark 3.1.} {\it From Lemma 3.1, if we take $r$
as the time $t$, then we know that the main part of $\mathcal {L}$
is like the seconder order operator $\p_t^2-\ds\f{1}{(1+t)^{2(\g-1)}}\Delta+\ds\f{2(\g-1)}{1+t}\p_t$,
which is strictly hyperbolic but degenerate as
$t\to\infty$. Recently, with respect to the semilinear wave equations with the forms
of $\p_t^2u-\Delta u+\ds\f{\mu}{(1+t)^{\al}}\p_tu=f(u)$, where $\mu>0$ and $\al>0$ are suitable constants, there
have been extensive
and interesting works on the global existence or blowup results for the different nonlinear function $f(u)$, one can see
[8], [23-24] and the references therein.}

\vskip 0.5 true cm
\centerline{\bf $\S 4$. A first-order weighted energy estimate}
\vskip 0.5 true cm

In this section, we establish a weighted energy estimate of
$\nabla_x\dot\Phi$ for the linear part of (3.1)
together with (3.4)-(3.5),
which will play a fundamental role in our subsequent analysis.

Set $D_T=\{(r,\th,\vp): 1<r<T, 0\le\th\le 2\pi, 0\le\vp<\vp_0\}$
for any $T>1$,  $B_T=\{(r,\th,\vp)\in \Sigma: 1<r<T,
0\le\th\le2\pi, \vp=\vp_0\}$, and $S_T=\bar\O\cap\{r=T\}$.

{\bf Theorem 4.1.} {\it Let $\dot\Phi\in C^{2}(\bar D_T)$ satisfy the boundary
condition (3.4) and initial data condition (3.5). Then there exists a multiplier
$\mathcal {M}\dot\Phi=r^{\mu}a(r)\p_r\dot\Phi$ such that for fixed
constant $\mu=4\g-6$ we have
\begin{align*}
&T^{\mu}\int_{S_T}(\p_r\dot\Phi)^2dS+T^{\mu-2\g}\int_{S_T}(Z\dot\Phi)^2dS
+C\int_{D_T}\big(r^{\mu-1-\dl}(\p_r\dot\Phi)^2+r^{\mu-1-2\g}(Z\dot\Phi)^2\big)dx\\
&\le \int_{D_T}\mathcal {L}\dot\Phi\cdot\mathcal {M}\dot\Phi
dx
+C\ve^2,\tag{4.1}
\end{align*}
where $(Z\dot\Phi)^2=\ss_{k=1}^3(Z_k\dot\Phi)^2$, $C>0$ is a generic positive constant, and $\dl>0$ is a fixed constant.}

{\bf Remark 4.1.} {\it Here we emphasize that the choice of $\mu=4\g-6$ in (4.1) is very necessary due to the following two reasons:
First, to guarantee the positivity of (4.3) below, one should let $\mu\le 4\g-6$; Second, by the Bernoulli's law (1.2),
we have $c^2(\rho)=
c^2(\hat\rho)-\f{\g-1}{2}\big((\p_r\dot\Phi)^2+2\hat U\p_r\dot\Phi+\f{1}{r^2}(Z\dot\Phi)^2\big)$.
Notice that only the estimate of $|\nabla_x\dot\Phi|\le
C\ve r^{-\f{\mu}{2}-1}$ can be obtained by the analysis in $\S 5-\S 6$. On the other hand, $c^2(\hat\rho(r))\ge Cr^{-2(\g-1)}$ and $\hat U=O(1)$ hold
by Lemma 2.1. Therefore, in order to guarantee the absence of vacuum for any finite $r$ in $\O$, we require to choose
the constant $\mu$ such that $-\ds\f{\mu}{2}-1\le
-2(\g-1)$, which leads to $\mu\ge 4\g-6.$ Combining these two reasons yields $\mu=4\g-6$.}

{\bf Remark 4.2.} {\it In Theorem 4.1, it suffices to choose the constant $\dl>0$. However, to derive the higher order energy estimates
of $\dot\Phi$, we require to give more restrictions on $\dl$ (one can see Theorem 5.1 in $\S 5$).}

{\bf Proof.} It follows from the integration by parts and (3.4)-(3.5) that
\begin{align*}
&\int_{D_T}\mathcal {L}\dot\Phi\cdot\mathcal {M}\dot\Phi
dx\\
&=
\int_{S_T}\f{1}{2}r^{\mu}a(r)(\p_r\dot\Phi)^2dS+\int_{S_T}\f{1}{2}r^{\mu-2}
P_1a(r)(Z\dot\Phi)^2dS\\
&\quad -\bigg(\int_{S_1}\f{1}{2}r^{\mu}a(r)(\p_r\dot\Phi)^2dS+\int_{S_1}\f{1}{2}r^{\mu-2}
P_1a(r)(Z\dot\Phi)^2dS\bigg)\\
&\quad +\int_{D_T}\bigg(r^{\mu-1}\big((P_2-\f{\mu+2}{2})a(r)-\f{1}{2}ra'(r)\big)(\p_r\dot\Phi)^2
-\f{1}{2}r^{\mu-3}\big((\mu P_1+rP'_1)a(r)+ra'(r)P_1\big)
(Z\dot\Phi)^2\bigg)dx\\
&\quad +\int_{B_T}r^{\mu}a(r)P_1(-x_2Z_2\dot\Phi+x_1Z_3\dot\Phi)\p_r\dot\Phi\\
&\ge \int_{S_T}\f{1}{2}r^{\mu}a(r)(\p_r\dot\Phi)^2dS+\int_{S_T}\f{1}{2}r^{\mu-2}
P_1a(r)(Z\dot\Phi)^2dS\\
&\quad +\int_{D_T}\bigg(r^{\mu-1}\big((P_2-\f{\mu+2}{2})a(r)-\f{1}{2}ra'(r)\big)(\p_r\dot\Phi)^2
-\f{1}{2}r^{\mu-3}\big((\mu P_1+rP'_1)a(r)+ra'(r)P_1\big)
(Z\dot\Phi)^2\bigg)dx\\
&\quad -C\ve^2\tag{4.2}
\end{align*}

It is noted that
\begin{align*}
&(P_2-\f{\mu+2}{2})a(r)-\f{1}{2}ra'(r)\\
&=\f{1}{2(\hat{U}^2-{c}^2(\hat\rho))}\big((4\g-6-\mu)\hat{U}^4+(4+2(\mu+2))\hat{U}^2{c}^2(\hat\rho)
+(2-\mu){c}^4((\hat\rho))\big)a(r)-\f{1}{2}ra'(r),\tag{4.3}
\end{align*}
then in order to guarantee the positivity of (4.3) for $\mu=4\g-6$, we require
$$a(r)>0\quad \text{and $a'(r)<0$}.\eqno{(4.4)}$$

For this end, we choose
$$a(r)=1+r^{-\dl}\quad\text{with $\dl>0$.}$$

In this case, one can arrive at
$$
(P_2-\f{\mu+2}{2})a(r)-\f{1}{2}ra'(r)>\f{1}{2}\dl r^{-\dl}.\eqno{(4.5)}
$$

On the other hand, it follows from a direct computation and the assumption of $1<\g<2$ that
\begin{align*}
&-(\mu P_1+r\p_rP_1)a(r)-ra'(r)P_1\\
&=\f{c^2(\hat\rho)}{({\hat U}^2-{c}^2(\hat\rho))^3}\biggl((2(\g-1)-\mu){\hat U}^4
+(4+2\mu){\hat U}^2{c}^2(\hat\rho)
-\mu{c}^4(\hat\rho)\biggr)a(r)+\f{\dl {c}^2(\hat\rho)r^{-\dl}}{{\hat U}^2-{c}^2(\hat\rho)}\\
&>\f{2{c}^2(\hat\rho)}{({\hat U}^2-{c}^2(\hat\rho))^3}(2-\g){\hat U}^4\\
&>Cr^{-2(\g-1)}.\tag{4.6}
\end{align*}
Thus, substituting (4.5)-(4.6) into (4.2) yields Theorem 4.1. \qquad \qquad \qquad \qquad \qquad
 $\square$

\vskip 0.5 true cm
\centerline{\bf $\S 5$. Higher-order weighted energy estimates of $\dot\Phi$}
\vskip 0.5 true cm

In this section, we will derive the higher-order energy
estimates of solution $\dot\Phi$ to (3.1) with (3.4)-(3.5) so that the suitable decay properties of
$\nabla_x\dot\Phi$ can be obtained and the density $\rho(x)>0$ can be also derived in subsequent $\S 6$.
Due to the Neumann boundary condition (3.4),
the asymptotic degeneracy of some coefficients in (3.1), and the different decay rates of $\p_r\dP$ and
$\ds\f{Z\dP}{r}$, the related derivation procedure will become rather complicated
and technical.

{\bf Theorem 5.1.} {\it Let $\dot\Phi\in C^4(\bar{D}_T)$ be the solution to
(3.1) with (3.4)-(3.5), and further assume
$$
\sum\limits_{0\le l_1+l_2\le 1}
r^{l_1}|\p_r^{l_1}Z^{l_2}\p_r\dot\Phi|\le M\ve r^{-2(\g-1)},\quad
r^{-1}|Z\dot\Phi|\le M\ve r^{-\si},\quad
r^{-1}|Z^2\dot\Phi|\le M\ve r^{-(\g-1)},\eqno{(5.1)}
$$
where $M>0$ is a constant, and $\si=\min\{1, 2(\g-1)\}$. Then for sufficiently small $\ve>0$ and
$0\le k\le 3$, we have
\begin{align*}
&T^{\mu+2k}\int_{S_T}|\nabla_x^k\p_r\dot\Phi|^2dS
+T^{\mu-2\g+2k}\int_{S_T}|\nabla_x^k
(\f{1}{r}Z\dot\Phi)|^2dS\\
&\quad +\int_{D_T}\biggl(r^{\mu-1-\dl+2k}|\nabla_x^k\p_r\dot\Phi|^2+r^{\mu+1-2\g+2k}
|\nabla_x^k(\f{1}{r}Z\dot\Phi)|^2\biggr)dx\\
&\le
C\ve^2,\tag{5.2}
\end{align*}
where $\mu=4\g-6$, $0<\dl\le\min\{\g-1, \si-(\g-1)\}$, and the domains $D_T, B_T, S_T$ have been defined
in the beginning of $\S 4$.}

In order to prove Theorem 5.1, we will apply the induction method on $k$ in (5.2) to establish
the following estimates respectively:

(i) $\p_rS^k\dP$ and $ZS^k\dP$ with $S=r\p_r$ and $1\le k\le 3$ (in this case, all the radial derivatives of $\na_x\dP$
up to third order are treated);

(ii) $\p_rZ\dP$ and $Z^2\dP$ (in this case, together with the case $k=1$ in (i),
all the second order derivatives $\na_x^2\dot\Phi$ are treated);

(iii) $\p_rSZ\dP$, $ZSZ\dP$, $\p_rZ^2\dP$ and $Z^3\dP$ (in this case, together with the case $k=2$ in (i),
all the third order derivatives $\na_x^3\dot\Phi$ are treated);

(iv) $\p_rS^2Z\dP$, $ZS^2Z\dP$, $\p_rSZ^2\dP$, $ZSZ^2\dP$, $\p_rZ^3\dP$ and $Z^4\dP$ (in this case,
together with the case $k=3$ in (i),
all the fourth order derivatives $\na_x^4\dot\Phi$ are treated).

These estimates will be given in Lemma 5.2-Lemma 5.5 respectively.

At first, we establish the radial derivative estimates of $\dot\Phi$ under the suitable induction
assumption. Set $S=r\p_r$, which is tangent to
fixed nozzle wall $\Sigma$, then we have

{\bf Lemma 5.2. (Radial derivative estimates)} {Under the assumptions of Theorem 5.1, if (5.2)
holds for $0\le l\le m-1$ with $1\le m\le 3$,
then
\begin{align*}
&T^{\mu}\int_{S_T}(\p_r S^m\dot\Phi)^2dS+T^{\mu-2\g}\int_{S_T}(ZS^m\dot\Phi)^2dS
+\int_{D_T}\bigg(r^{\mu-1-\dl}(\p_r S^m\dot\Phi)^2+r^{\mu-1-2\g}(Z S^m\dot\Phi)^2\bigg)dx\\
&\le C\ve^2+ C\ve\int_{D_T}
\ds\sum_{l=0}^m\bigg(r^{\mu+1-2\g+2l}(\nabla_x^l(\f{1}{r}Z\dot\Phi))^2+r^{\mu-1-\dl+2l}(\nabla_x^l\p_r\dot\Phi)^2\bigg)
dx\\
&\quad +C\ve\ss_{l=0}^m\bigg(T^{\mu+2l}\int_{S_T}(\nabla_x^l\p_r\dot\Phi)^2dS+T^{\mu-2\g+2l}\int_{S_T}
(\nabla_x^l(Z\dot\Phi))^2dS\bigg)\\
&\quad +C\ve\bigg(\int_{D_T}r^{\mu-1-\dl}(\p_r S^m\dot\Phi)^2+r^{\mu-1-2\g}(Z S^m\dot\Phi)^2dx\bigg)^{\f12},\tag{5.3}
\end{align*}
where $0<\dl\le\g-1$.

Especially, for $m=0$, the following estimate holds
\begin{align*}
&T^{\mu}\int_{S_T}(\p_r\dot\Phi)^2dS+T^{\mu-2\g}\int_{S_T}(Z\dot\Phi)^2dS
+\int_{D_T}\bigg(r^{\mu-1-\dl}(\p_r \dot\Phi)^2+r^{\mu-1-2\g}(Z \dot\Phi)^2\bigg)dx\\
&\le C\ve^2+ C\ve\int_{D_T}
\bigg(r^{\mu+1-2\g}(\f{1}{r}Z\dot\Phi)^2+r^{\mu-1-\dl}(\p_r\dot\Phi)^2\bigg)
dx\\
&\quad +C\ve\bigg(T^{\mu}\int_{S_T}(\p_r\dot\Phi)^2dS+T^{\mu-2\g}\int_{S_T}
(Z\dot\Phi)^2dS\bigg).\tag{5.4}
\end{align*}}

{\bf Remark 5.1.} {\it For the case of $m=0$ in (5.4), we do not require any induction assumption.}

{\bf Remark 5.2.} {\it It is noted that the angular derivatives of $\dP$
are still included in the right hand side of (5.3), which implies that we have not
obtained the complete estimates on the radial derivative estimates of $\dP$. However,
since the coefficients of angular derivatives of $\dP$  in (5.3) are small, then together
with the subsequent angular derivative estimates, we can derive (5.2).}

{\bf Proof.} Noticing that on $\Sigma$
$$S^m Z_2\dot\Phi=S^m Z_3\dot\Phi=0.\eqno{(5.5)}
$$

This, together with Theorem 4.1 and (3.5), yields
\begin{align*}
&T^{\mu}\int_{S_T}(\p_rS^m\dot\Phi)^2dS+T^{\mu-2\g}\int_{S_T}(ZS^m\dot\Phi)^2dS
+\int_{D_T}\biggl(r^{\mu-1-\dl}(\p_rS^m\dot\Phi)^2+r^{\mu-1-2\g}(ZS^m\dot\Phi)^2\biggr)dx\\
&\le C\int_{D_T}\mathcal {L}S^m\dot\Phi\cdot\mathcal {M}S^m\dot\Phi
dx+C\ve^2.\tag{5.6}
\end{align*}

Next, we derive an explicit representation of
$\mathcal {L}S^m\dot\Phi$ for the later uses.

By a direct computation, we have
$$
\mathcal {L} S=S\mathcal {L}+2\mathcal {L}+\ds\f{P'_1}{r}\sum\limits_{i=1}^3Z_i^2+A_1,
\eqno{(5.7)}
$$
where $A_1=-P'_2\p_r$ is a first order operator.

By induction, for $1\le m\le3$, we further arrive at
$$
\mathcal {L}S^m=S^m\mathcal {L}+mS^{m-1}\bigl(\ds\f{P'_1}{r}\sum\limits_{i=1}^3Z_i^2\big)
+\sum\limits_{0\le l\le m-1}C_{lm}S^{l}\mathcal{L}\dot\Phi+A_m,\eqno{(5.8)}
$$
where $C_{lm}$ are some suitable constants, $A_m$ stands for a lower order differential operator whose order is less
than $m$. For examples,
\begin{align*}
&A_2=\big(\f{2P'_1}{r}-r(\f{P'_1}{r})'\big)\sum\limits_{i=1}^3Z_i^2
+SA_1+A_1S+A_1,\\
&A_3=SA_2+[\f{P'_1}{r}\sum\limits_{i=1}^3Z_i^2, S^2]
+2S(\f{P'_1}{r}\sum\limits_{i=1}^3Z_i^2)+A_2+A_1S^2,
\end{align*}
here and below, $[\cdot , \cdot]$ denotes the usual commutator.

For convenient treatments, for $1\le m\le3$, we rewrite (5.8) as
$$\mathcal {L}S^m=S^m\mathcal {L}+B_{1m}+B_{2m}\eqno{(5.9)}$$
with
\begin{align*}
&B_{1m}=\ds\sum_{0\le l\le
m-1}C_{lm} S^{l}\mathcal{L}+ \ds\sum_{0\le l\le
m-2}C_{lm}S^{m-1-l}\big(\f{rP'_1}{P_1}\big)S^l\big(\f{P_1}{r^2}\ss_{i=1}^3Z_i^2\big)+A_m,\\
&B_{2m}=\f{mrP'_1}{P_1}S^{m-1}\big(\f{P_1}{r^2}\ds\sum_{i=1}^3Z_i^2\big),
\end{align*}
where $B_{2m}\dP$ contains the $(m+1)-$th order (the highest order) derivatives of $\dot\Phi$,
but $B_{1m}\dP$ only includes $\na_x^{\al}\dP$ with $|\al|\le m$ (the lower order derivatives of $\dP$)
and  $\na_x^{m+1}\dP$ with small coefficients.

In addition, from the equation (3.1), for $0\le m\le3$, we have
$$
S^m\mathcal{L}\dot\Phi=I_{1}^m+I_{2}^m+I_{3}^m,\eqno{(5.10)}
$$
where
\begin{align*}
&I_{1}^m=f_{00}\p_r^2S^m\dot\Phi+\f{1}{r^2}\sum\limits_{1\le
i,j\le 3}f_{ij}Z_iZ_jS^m\dot\Phi
+\f{1}{r}\sum\limits_{i=1}^3f_{0i}\p_rZ_iS^m\dot\Phi,\\
&I_{2}^m=f_{00}[S^m, \p_r^2]\dot\Phi
+\f{1}{r}\sum\limits_{i=1}^3f_{0i}[S^m, \p_rZ_i]\dot\Phi,\\
&I_{3}^m=\sum\limits_{0\le l\le
m}C_{lm}\big\{\sum\limits_{l_1+l_2=
l, l_1\ge1}\t C_{l_1l_2}\big(S^{l_1}f_{00}S^{l_2}\p_r^2\dot\Phi
+S^{l_1}(\f{1}{r^2}\sum\limits_{1\le i,j\le3}f_{ij})S^{l_2}Z_iZ_j\dot\Phi\\
&\qquad\quad +S^{l_1}(\f{1}{r}\sum\limits_{i=1}^3f_{0i})S^{l_2}\p_rZ_i\dot\Phi\big)\big\}+S^mf_0.\\
\end{align*}

Based on the preparations above, we now treat $\int_{D_T}\mathcal {L}S^m\dot\Phi\cdot\mathcal{M}S^m\dot\Phi$
in the right hand side of (5.6).
This procedure is divided into the following
five parts.

\vskip 0.3 true cm

{\bf Part 1. The estimate of $\int_{D_T}I_{1}^m\cdot\mathcal{M}S^m\Phi dx$}

\vskip 0.3 true cm

Notice that we have for $C^1-$smooth functions $g_i$ ($1\le i\le 3$)
$$
\sum\limits_{i=1}^3Z_ig_i
=\p_1(x_3g_3-x_2g_1)+\p_2(x_1g_1-x_3g_2)+\p_3(x_2g_2-x_1g_3)\eqno{(5.11)}
$$
and
$$
\int_{D_T}\sum\limits_{i=1}^3Z_ig_idx=-\int_{B_T}\f{1}{\sin\vp_0}(x_2g_2-x_1g_3)dS.\eqno{(5.12)}
$$

In addition, a direct computation yields for $m\le3$
\begin{align*}
I_{1}^m&\cdot \mathcal{M}S^m\dot\Phi\\
=&\p_r\bigg(\f{1}{2}r^{\mu}a(r)f_{00}(\p_rS^m\dP)^2-r^{\mu-2}a(r)\ds\sum_{1\le i<j\le3}f_{ij}Z_iS^m\dP Z_jS^m\dP
-\f{1}{2}r^{\mu-2}a(r)\ds\sum_{i=1}^3f_{ii}(Z_iS^m\dP)^2\bigg)\\
&+\ds\sum_{i=1}^3Z_i\bigg(\f{1}{2}r^{\mu-1}a(r)f_{0i}(\p_rS^m\dP)^2+r^{\mu-2}a(r)\p_rS^m\dP\ss_{j=1}^3f_{ij}Z_jS^m\dP\bigg)\\
&-\f{1}{2}\p_r(r^\mu a(r)f_{00})(\p_rS^m\dP)^2-\f{1}{2}r^{\mu-1}a(r)(\p_rS^m\dP)^2\ds\sum_{i=1}^3Z_if_{0i}
+\ds\sum_{i=1}^3\p_r(\f{1}{2}r^{\mu-2}a(r)f_{ii})(Z_iS^m\dP)^2\\
&+\ds\sum_{1\le i<j\le3}\p_r(r^{\mu-2}a(r)f_{ij})Z_iS^m\dP Z_jS^m\dP.\tag{5.13}
\end{align*}

On the other hand, by the expressions of $f_{ij}, f_{0i}$ and (5.5), a crucial observation yields on $B_T$
\begin{align*}
&x_2\bigg(\f{1}{2}r^{\mu-1}a(r)f_{02}(\p_rS^m\dP)^2
+r^{\mu-2}a(r)\p_rS^m\dP\ds\sum_{j=1}^3f_{2j}Z_jS^m\dP\bigg)\\
&\quad -x_1\bigg(\f{1}{2}r^{\mu-1}a(r)f_{03}(\p_rS^m\dP)^2
+r^{\mu-2}a(r)\p_rS^m\dP\ds\sum_{j=1}^3f_{3j}Z_jS^m\dP\bigg)\\
&=\f{1}{2}a(r)(\p_rS^m\dP)^2(x_2f_{02}-x_1f_{03})
+r^{\mu-2}a(r)\p_rS^m\dP\ds\sum_{j=1}^3(x_2f_{2j}-x_1f_{3j})Z_jS^m\dP\\
&=\f{1}{\hat{U}^2-{c}^2(\hat\rho)}r^{\mu-2}a(r)\p_rS^m\dP\bigg(\f{\g-1}{2}(\p_r\dP)^2
-(\g-1)\hat{U}\p_r\dP+\f{\g-1}{2r^2}\ds\sum_{k=1}^3(Z_k\dP)^2\bigg)
(x_2Z_2S^m\dP\\
&\quad -x_1Z_3S^m\dP)-\f{1}{\hat{U}^2-{c}^2(\hat\rho)}r^{\mu-4}a(r)\p_rS^m\dP\ss_{j=1}^3Z_j\dP Z_jS^m\dP(x_2Z_2\dP-x_1Z_3\dP)\\
&\quad -\f{1}{2}r^{\mu-2}a(r)(\hat{U}+\p_r\dP)(\p_rS^m\dP)^2(x_2Z_2\dP-x_1Z_3\dP)\\
&=0.\tag{5.14}
\end{align*}

Thus, by (5.13) together with (5.12) and (5.14), it follows from an integration by parts
and simultaneously notices
the expressions of $f_i$ and the assumption (5.1) that
\begin{align*}
&|\int_{D_T}I_{1}^m\cdot\mathcal{M}S^m\dP dx|\\
&\le C\ve^2+C\ve\bigg(T^{\mu}\int_{S_T}(\p_rS^m\dot\Phi)^2dS
+T^{\mu-2\g}\int_{S_T}(ZS^m\dot\Phi)^2dS\\
&\quad+\int_{D_T}\big(r^{\mu-1-\dl}(\p_rS^m\dot\Phi)^2
+r^{\mu-1-2\g}(ZS^m\dot\Phi)^2\big)dx\bigg),\tag{5.15}
\end{align*}
here we have used some facts such as
$$|r^{\mu-1}a(r)(\p_rS^m\dP)^2\ds\sum_{i=1}^3Z_if_{0i}|\le C\ve r^{\mu-1-(\g-1)}(\p_rS^m\dP)^2\le C\ve r^{\mu-1-\dl}(\p_rS^m\dP)^2\quad \text{for $0<\dl\le \g-1$}.$$

\vskip 0.3 true cm

{\bf Part 2. The estimate of $\int_{D_T}I_{2}^m\cdot\mathcal{M}S^m\Phi dx$}

\vskip 0.3 true cm

It follows from the expressions
of $f_i$, Lemma 2.5 (ii) and (5.1) that
$$|I_{2}^m|\le C\ve
\big(r^{-\si-1}\ss_{0\le l\le
m-1}|S^l\p_rZ\dP|+r^{-2(\g-1)}\ss_{0\le l\le
m-1}|S^l\p_r^2\dP|\big),$$
which derives that
$$\int_{D_T}|I_{2}^m\cdot\mathcal{M}S^m\dP|dx
\le C\ve\int_{D_T}\bigg(\ds\sum_{0\le l\le
m-1}r^{\mu-1-\dl}(S^l\p_rZ\dot\Phi)^2 +\ds\sum_{1\le l\le
m}r^{\mu-1-\dl}(S^l\p_r\dP)^2\bigg)dx.\eqno{(5.16)}
$$

\vskip 0.3 true cm

{\bf Part 3. The estimate of $\int_{D_T}I_{3}^m\cdot\mathcal{M}S^m\Phi dx$}

\vskip 0.3 true cm

At first, we treat the case of $\int_{D_T}I_{3}^m\cdot\mathcal{M}S^m\Phi dx$ with $m\le 2$.

For $m\le 2$, as in Part 2 it follows from the expressions
of $f_i$ and the assumption (5.1) that
$$ |I_{3}^m|\le C\ve \bigg(\ds\sum_{0\le
l\le m}\big(r^{-\g}|S^l\p_r\dP|+r^{-\si-2}|S^lZ\dP|\big)+\ds\sum_{0\le l\le m-1}
r^{-2\g}|S^lZ^2\dP|\bigg),$$
which derives for $m\le 2$
$$
\int_{D_T}|I_{3}^m\cdot\mathcal{M}S^m\dP|dx
\le C\ve\int_{D_T}
\ds\sum_{l=0}^m\bigg(r^{\mu+1-2\g+2l}(\nabla_x^l(\f{1}{r}Z\dP))^2+r^{\mu-1-\dl+2l}(\nabla_x^l\p_r\dP)^2\bigg)
dx.\eqno{(5.17)}$$

Next we deal with $\int_{D_T}I_{3}^3\cdot\mathcal{M}S^3\dP dx$.

It is noted that the most troublesome terms in $I_{3}^3$ are the ones which include the products of third order derivatives
of $\dot\Phi$ since there are no related weighted $L^{\infty}$ estimates in (5.1). For the convenient  treatments, we decompose $I_{3}^3$ into $J_1$ and $J_2$
by using $S^2=r\p_r+r^2\p_r^2$, where only $J_2$ contains the product terms of third order derivatives
of $\dot\Phi$. Namely,
$$I_{3}^3=J_1+J_2\eqno{(5.18)}$$
with
\begin{align*}
&J_1=\\
&\sum\limits_{0\le l\le
2}C_{l2}\sum\limits_{l_1+l_2=
l, l_1\ge1}\t C_{l_1l_2}\bigg(S^{l_1}f_{00}S^{l_2}\p_r^2\dot\Phi
+S^{l_1}(\f{1}{r^2}\sum\limits_{1\le i,j\le3}f_{ij})S^{l_2}Z_iZ_j\dot\Phi+S^{l_1}(\f{1}{r}\sum\limits_{i=1}^3f_{0i})S^{l_2}\p_rZ_i\dot\Phi\bigg)\\
&\quad +C_{33}\sum\limits_{(l_1,l_2)\not=(2, 1)}\t C_{l_1l_2}\bigg(S^{l_1}f_{00}S^{l_2}\p_r^2\dot\Phi
+S^{l_1}(\f{1}{r^2}\sum\limits_{1\le i,j\le3}f_{ij})S^{l_2}Z_iZ_j\dot\Phi+S^{l_1}(\f{1}{r}\sum\limits_{i=1}^3f_{0i})S^{l_2}\p_rZ_i\dot\Phi\bigg)\\
&\quad +C_{33}\t C_{21}\bigg(r\p_rf_{00}S\p_r^2\dP+r\p_r(\f{1}{r^2}\ss_{i,j=1}^3f_{ij})SZ_iZ_j\dP
+r\p_r(\f{1}{r}\ss_{i=1}^3f_{0i})S\p_rZ_i\dP\bigg)+S^3f_0
\end{align*}
and
$$J_2=C_{33}\t C_{21}\bigg(r^2\p_r^2f_{00}S\p_r^2\dP+r^2\p_r^2(\f{1}{r^2}\ds\sum_{i,j=1}^3f_{ij})SZ_iZ_j\dP
+r^2\p_r^2(\f{1}{r}\ds\sum_{i=1}^3f_{0i})S\p_rZ_i\dP\bigg).$$

By the assumption (5.1) and the expressions of $f_{ij}, f_0$, then a direct computation yields
$$|J_1|\le C\ve \bigg\{\ds\sum_{0\le
l\le 3}\big(r^{-\g}|S^l\p_r\dP|+r^{-\si-2}|S^lZ\dP|\big)+\ds\sum_{0\le l\le2}r^{-2\g}|S^lZ^2\dP|\biggr\}.\eqno{(5.19)}$$

Next, by the expressions of $f_{ij}$ we continue to decompose $J_2$  as $J_2=J_{21}+J_{22}$ so that only
$J_{22}$ contains the product terms of third order derivatives
of $\dot\Phi$. More concretely,
$$J_2=J_{21}+J_{22}$$
with
\begin{align*}
J_{22}&=-\f{r^2}{{\hat U}^2-c^2(\hat\rho)}\bigg\{\bigg((\g+1)\big(\hat {U}+\p_r\dP\big)\p_r^3\dP+\f{\g-1}{r^2}
\ds\sum_{i=1}^3Z_i\dP\p_r^2Z_i\dP\bigg)S\p_r^2\dP\\
&-\f{\g-1}{r^2}\ds\sum_{i=1}^3\bigg(\p_r\dP\p_r^3\dP-\hat{U}\p_r^3\dP+\f{1}{r^2}\ds\sum_{k=1}^3Z_k\dP\p_r^2Z_k\dP
-\f{1}{r^2}Z_i\dP\p_r^2Z_i\dP\bigg)SZ_i^2\dP\\
&+\f{1}{r^4}\ds\sum_{1\le i\not=j\le 3}\big(Z_i\dP\p_r^2Z_j\dP+Z_j\dP\p_r^2Z_i\dP\big)SZ_iZ_j\dP +\f{1}{r^2}\ds\sum_{i=1}^3\big(\hat{U}\p_r^2Z_i\dP+\p_r^3\dP Z_i\dP\\
&+\p_r\dP\p_r^2Z_i\dP\big)S\p_rZ_i\dP\biggr\}
\end{align*}
and
$$|J_{21}|\le C\ve \bigg(r^{-(\g-1)}|S\p_r^2\dP|+r^{-\si-1}|S\p_rZ\dP|+r^{-2\g}|SZ^2\dP|\bigg),\eqno{(5.20)}$$
here we point out that (5.20) is derived by  direct but tedious computations through applying Lemma 2.1, assumption
(5.1) and the concrete expression of $J_{21}$.

Combining (5.19) and (5.20) together with Lemma 2.5 can yield
$$
\int_{D_T}|(J_1+J_{21})\cdot\mathcal{M}S^3\dP|dx
\le C\ve\int_{D_T}
\ds\sum_{l=0}^3\bigg(r^{\mu+1-2\g+2l}|\nabla_x^l(\f{1}{r}Z\dP)|^2
+r^{\mu-1-\dl+2l}|\nabla_x^l\p_r\dP|^2\bigg)dx.\eqno{(5.21)}
$$

Finally we treat $\int_{D_T}|J_{22}\cdot\mathcal {M}S^3\dP|dx$.

To overcome the difficulties induced by the lack of weighted $L^{\infty}$ estimates of
$|\na_x^3 \dot\Phi|$ in $J_{22}$, we will use the interpolation inequalities in Corollary 2.3 and Lemma 2.6.
In fact, by (5.1) and the expression of $J_{22}$, it is only enough to deal with the following typical terms
in $\int_{D_T}|J_{22}\cdot\mathcal {M}S^3\dP|dx$:
\vskip 0.3 true cm
{\bf (A) Estimate of $|r^{\mu+2}\p_r^3\dP S\p_r^2\dP\p_rS^3\dP|_{L^1(D_T)}$}
\begin{align*}
&|r^{\mu+2}\p_r^3\dP S\p_r^2\dP\p_rS^3\dP|_{L^1}=|r^{\mu+3}(\p_r^3\dP)^2\p_rS^3\dP|_{L^1}\\
=&|r^{\dl-2(\g-1)}\cdot (r^{\f{8\g-3-\dl}{4}}\p_r^3\dP)^2\cdot r^{\f{\mu-1-\dl}{2}}\p_rS^3\dP|_{L^1}\\
\le & |r^{\f{8\g-3-\dl}{4}}\p_r^3\dP|_{L^4}^2|r^{\f{\mu-1-\dl}{2}}\p_rS^3\dP|_{L^2}\\
\le & C\ve\ds\sum_{l=0}^3\bigg(|r^{\f{\mu-1-\dl+2l}{2}}\nabla_x^l\p_r\dP|_{L^2}
+|r^{\f{\mu+1-2\g+2l}{2}}\nabla_x^l(\f{1}{r}Z\dP)|_{L^2}\bigg)
|r^{\f{\mu-1-\dl}{2}}\p_rS^3\dP|_{L^2}\\
&\qquad \quad\qquad \quad\qquad \quad\qquad \quad\qquad \quad\text{(Applying Lemma 2.7 for $\dP$)}\\
\le & C\ve\ds\sum_{l=0}^3\bigg(|r^{\f{\mu-1-\dl+2l}{2}}\nabla_x^l\p_r\dP|_{L^2}^2
+|r^{\f{\mu+1-2\g+2l}{2}}\nabla_x^l(\f{1}{r}Z\dP)|_{L^2}^2\bigg).\tag{5.22}
\end{align*}

\vskip 0.3 true cm
{\bf (B) Estimate of $|r^{\mu}Z\dP\p_r^2Z\dP S\p_r^2\dP\p_rS^3\dP|_{L^1(D_T)}
+|r^{\mu}Z\dP\p_r^3\dP S\p_rZ\dP\p_rS^3\dP|_{L^1(D_T)}$}
\begin{align*}
&|r^{\mu}Z\dP\p_r^2Z\dP S\p_r^2\dP\p_rS^3\dP|_{L^1(D_T)}+|r^{\mu}Z\dP\p_r^3\dP S\p_rZ\dP\p_rS^3\dP|_{L^1(D_T)}\\
=&2|r^{\mu+1}Z\dP\p_r^3\dP\p_r^2Z\dP\p_rS^3\dP|_{L^1(D_T)}\\
\le& C|r^{\mu+2-\si}\p_r^3\dP\p_r^2Z\dP\p_rS^3\dP|_{L^1(D_T)} \qquad\qquad \qquad\text{(By assumption (5.1))}\\
=&|r^{\dl-2(\g-1)-\si}\cdot r^{\f{8\g-3-\dl}{4}}\p_r^3\dP\cdot r^{\f{8\g-7-\dl}{4}}\p_r^2Z\dP
\cdot r^{\f{\mu-1-\dl}{2}}\p_rS^3\dP|_{L^1(D_T)}\\
\le& C|r^{\f{8\g-3-\dl}{4}}\p_r^3\dP|_{L^4(D_T)}|r^{\f{8\g-7-\dl}{4}}\p_r^2Z\dP|_{L^4(D_T)} |r^{\f{\mu-1-\dl}{2}}\p_rS^3\dP|_{L^2(D_T)}\\
\le&C\ve\ds\sum_{l=0}^3\bigg(|r^{\f{\mu-1-\dl+2l}{2}}\nabla_x^l\p_r\dP|_{L^2(D_T)}^2
+|r^{\f{\mu+1-2\g+2l}{2}}\nabla_x^l(\f{1}{r}Z\dP)|_{L^2(D_T)}^2\bigg).\quad\text{(Applying Lemma 2.7 for $\dP$)}\tag{5.23}
\end{align*}

\vskip 0.3 true cm
{\bf (C) Estimate of $|r^{\mu}\p_r^3\dP SZ^2\dP\p_rS^3\dP|_{L^1(D_T)}$}
\begin{align*}
&|r^{\mu}\p_r^3\dP SZ^2\dP\p_rS^3\dP|_{L^1(D_T)}
=|r^{\mu+1}\p_r^3\dP\p_rZ^2\dP\p_rS^3\dP|_{L^1(D_T)}\\
=&|r^{\dl-2(\g-1)}\cdot r^{\f{8\g-3-\dl}{4}}\p_r^3\dP\cdot r^{\f{8\g-11-\dl}{4}}\p_rZ^2\dP
\cdot r^{\f{\mu-1-\dl}{2}}\p_rS^3\dP|_{L^1(D_T)}\\
\le& C|r^{\f{8\g-3-\dl}{4}}\p_r^3\dP|_{L^4(D_T)}|r^{\f{8\g-11-\dl}{4}}\p_rZ^2\dP|_{L^4(D_T)} |r^{\f{\mu-1-\dl}{2}}\p_rS^3\dP|_{L^2(D_T)}\\
\le&C\ve\ds\sum_{l=0}^3\bigg(|r^{\f{\mu-1-\dl+2l}{2}}\nabla_x^l\p_r\dP|_{L^2(D_T)}^2
+|r^{\f{\mu+1-2\g+2l}{2}}\nabla_x^l(\f{1}{r}Z\dP)|_{L^2(D_T)}^2\bigg).\quad\text{(Applying Lemma 2.7 for $\dP$)}\tag{5.24}
\end{align*}

\vskip 0.3 true cm
{\bf (D) Estimate of $|r^{\mu-2}Z\dP \p_r^2Z\dP SZ^2\dP\p_rS^3\dP|_{L^1(D_T)}$}
\begin{align*}
&|r^{\mu-2}Z\dP \p_r^2Z\dP SZ^2\dP\p_rS^3\dP|_{L^1(D_T)}\\
\le &C|r^{\mu-\si}\p_r^2Z\dP\p_rZ^2\dP\p_rS^3\dP|_{L^1(D_T)}\qquad\qquad\qquad\qquad\text{(By assumption (5.1))}\\
=& C|r^{\dl-2(\g-1)-\si}\cdot r^{\f{8\g-7-\dl}{4}}\p_r^2Z\dP\cdot r^{\f{8\g-11-\dl}{4}}\p_rZ^2\dP
\cdot r^{\f{\mu-1-\dl}{2}}\p_rS^3\dP|_{L^1(D_T)}\\
\le& C|r^{\f{8\g-7-\dl}{4}}\p_r^2Z\dP|_{L^4(D_T)}|r^{\f{8\g-11-\dl}{4}}\p_rZ^2\dP|_{L^4(D_T)} |r^{\f{\mu-1-\dl}{2}}\p_rS^3\dP|_{L^2(D_T)}\\
\le& C\ve\ds\sum_{l=0}^3\bigg(|r^{\f{\mu-1-\dl+2l}{2}}\nabla_x^l\p_r\dP|_{L^2(D_T)}^2
+|r^{\f{\mu+1-2\g+2l}{2}}\nabla_x^l(\f{1}{r}Z\dP)|_{L^2(D_T)}^2\bigg).\quad\text{(Applying Lemma 2.7 for $\dP$)}\tag{5.25}
\end{align*}

\vskip 0.3 true cm
{\bf (E) Estimate of $|r^{\mu}\p_r^2Z\dP S\p_rZ\dP\p_rS^3\dP|_{L^1(D_T)}$}
\begin{align*}
&|r^{\mu}\p_r^2Z\dP S\p_rZ\dP\p_rS^3\dP|_{L^1(D_T)}\\
= &|r^{\mu+1}(\p_r^2Z\dP)^2\p_rS^3\dP|_{L^1(D_T)}\\
=& |r^{\dl-2(\g-1)}\cdot (r^{\f{8\g-7-\dl}{4}}\p_r^2Z\dP)^2\cdot r^{\f{\mu-1-\dl}{2}}\p_rS^3\dP|_{L^1(D_T)}\\
\le &C|r^{\f{8\g-7-\dl}{4}}\p_r^2Z\dP|_{L^4(D_T)}^2|r^{\f{\mu-1-\dl}{2}}\p_rS^3\dP|_{L^2(D_T)}\\
\le &C\ve\ds\sum_{l=0}^3\bigg(|r^{\f{\mu-1-\dl+2l}{2}}\nabla_x^l\p_r\dP|_{L^2(D_T)}^2
+|r^{\f{\mu+1-2\g+2l}{2}}\nabla_x^l(\f{1}{r}Z\dP)|_{L^2(D_T)}^2\bigg).\quad\text{(Applying Lemma 2.7 for $\dP$)}
\tag{5.26}
\end{align*}

Collecting (5.22)-(5.26), one has
$$\int_{D_T}|J_{22}\cdot\mathcal{M}S^3\dP|dx
\le C\ve\int_{D_T}
\ds\sum_{l=0}^3\bigg(r^{\mu+1-2\g+2l}(\nabla_x^l(\f{1}{r}Z\dP))^2+r^{\mu-1-\dl+2l}(\nabla_x^l\p_r\dP)^2\bigg)
dx.\eqno{(5.27)}$$

This, together with (5.17) and (5.21), yields for $m\le 3$
$$\int_{D_T}|I_{3}^m\cdot\mathcal{M}S^m\dP|dx
\le C\ve\int_{D_T}
\ds\sum_{l=0}^m\bigg(r^{\mu+1-2\g+2l}(\nabla_x^l(\f{1}{r}Z\dP))^2+r^{\mu-1-\dl+2l}(\nabla_x^l\p_r\dP)^2\bigg)
dx.\eqno{(5.28)}$$

\vskip 0.3 true cm

{\bf Part 4. The estimate of $\int_{D_T}B_{1m}\dot\Phi\cdot\mathcal{M}S^m\Phi dx$}
\vskip 0.3 true cm

At first, from the expressions of $A_k$ with $1\le k\le m$ in (5.8) and Lemma 2.1, we have
\begin{equation}
\left\{
\begin{aligned}
&|A_1\dP|\le Cr^{1-2\g}|\p_r\dP|,\\
&|A_k\dP|\le Cr^{1-2\g}\big(\ss_{l=0}^{k-1}|\p_rS^l\dP|+\ss_{l=0}^{k-2}|\f{1}{r}S^lZ^2\dP|\big),\qquad
k=2,3.
\end{aligned}
\right.\tag{5.29}
\end{equation}
Substituting (5.29) into $B_{1m}\dP$ and using the expression of $\dot f$, then we have by (5.1)
and a direct computation
\begin{align*}
&|B_{11}\dP|\le C\ve\bigg(\ss_{0\le l\le 1}(r^{-\g}|S^l\p_r\dP|+r^{-\si-2}|S^lZ\dP|)+r^{-2\g}|Z^2\dP|\bigg)
+Cr^{1-2\g}|\p_r\dP|,
\tag{5.30}\\
&|B_{1k}\dP|\le C\ve\bigg(\ss_{0\le l\le k}(r^{-\g}|S^l\p_r\dP|+r^{-\si-2}|S^lZ\dP|)
+\ss_{0\le l\le k-1}r^{-2\g}|S^lZ^2\dP|\bigg)\\
&\quad \qquad +Cr^{1-2\g}(\ss_{l=0}^{k-1}|\p_rS^l\dP|+\ss_{l=0}^{k-2}|\f{1}{r}S^lZ^2\dP|),\qquad k=2,3.\tag{5.31}
\end{align*}

Since (5.2) holds for $l\le m-1$, we then have from (5.30)-(5.31)
\begin{align*}
&\int_{D_T}|B_{1m}\dot\Phi\cdot\mathcal{M}S^m\Phi dx|\le
C\ve\int_{D_T}\ss_{l=0}^m\bigg(r^{\mu+1-2\g+2l}(\na_x(\f{Z\dP}{r}))^2+r^{\mu-1-\dl+2l}
(\na_x^l\p_r\dP)^2\bigg)dx\\
&\qquad +C\ve\bigg(\int_{D_T}r^{\mu-1-\dl}(\p_rS^m\dP)^2dx\bigg)^{\f12}.\tag{5.32}
\end{align*}

\vskip 0.3 true cm

{\bf Part 5. The estimate of $\int_{D_T}B_{2m}\dot\Phi\cdot\mathcal{M}S^m\Phi dx$}
\vskip 0.3 true cm

It is noted that $B_{2m}\dot\Phi$ contains the $(m+1)-$th (the highest order) derivatives of $\dP$
and then $B_{2m}\dot\Phi\cdot\mathcal{M}S^m\Phi$ will contain the term $\na_{x}^{\al}\dP\na_{x}^{\beta}\dP$
($|\al|=|\beta|=m+1$) which will yield serious troubles in the general case.
However, thanks to $\ds\f{P_1'(r)}{P_1(r)}<0$
given in Lemma 3.1 and the good form of (3.1), the bad influence of $\na_{x}^{\al}\dP\na_{x}^{\beta}\dP$
with $|\al|=|\beta|=m+1$  can be eliminated in the related energy
estimates. We now give the details.

Since
$$\p_r S^m\dP=\ds\sum_{0\le l\le
m-1}C_{lm}\p_r S^l\dP+r S^{m-1}\p_r^2\dP,\eqno{(5.33)}$$
then it
follows from (3.1) that
$$\p_r S^m\dP=\ds\sum_{0\le l\le
m-1}C_{lm}\p_r S^l\dP+r S^{m-1}\big(\f{1}{r^2}P_1\ss_{i=1}^3Z_i^2\dP+\mathcal{L}\dP
-\f{1}{r}P_2\p_r\dP\big).\eqno{(5.34)}$$

A direct computation yields
\begin{align*}
&\int_{D_T} B_{2m}\dP\cdot \mathcal{M} S^m\dP dx=\int_{D_T}
\f{mrP'_1}{P_1}S^{m-1}\big(\f{P_1}{r^2}\sum\limits_{i=1}^3Z_i^2\dP\big)\cdot
r^{\mu}a(r)\p_r S^m\dP dx\\
= &\int_{D_T}
\f{mr^{\mu+1}a(r)P'_1}{P_1}\bigg\{\bigg[S^{m-1}\big(\f{P_1}{r^2}\sum\limits_{i=1}^3Z_i^2\dP\big)\bigg]^2\\
&\qquad + \bigg(\ds\sum_{0\le l\le
m-1}C_{lm}\p_r S^l\dP+r S^{m-1}\big(\mathcal{L}\dP-\f{P_2}{r}\p_r\dP\big)\bigg)
S^{m-1}\big(\f{P_1}{r^2}\sum\limits_{i=1}^3Z_i^2\dP\big)\bigg\}dx\\
&\le \int_{D_T} \f{mr^{\mu+1}a(r)P'_1}{P_1}\bigg(\ds\sum_{0\le l\le
m-1}C_{lm}\p_r S^l\dP+r S^{m-1}\big(\mathcal{L}\dP-\f{P_2}{r}\p_r\dP\big)\bigg)
S^{m-1}\big(\f{P_1}{r^2}\sum\limits_{i=1}^3Z_i^2\dP\big)dx.\\
&\quad\quad\quad\quad\quad\quad\quad\quad\quad\quad\quad\quad
\quad\quad\text{(By $\ds\f{P_1'(r)}{P_1(r)}<0$)}\tag{5.35}
\end{align*}

Note that $\ds\sum_{0\le l\le
m-1}C_{lm}\p_r S^l\dP-r S^{m-1}\big(\f{P_2}{r}\p_r\dP\big)$ only contains at most
$m-$order derivatives of $\dP$, then we have by (5.2) for $l\le m-1$
$$
\int_{D_T}r^{\mu-1-\dl}\bigg(\ds\sum_{0\le l\le
m-1}C_{lm}\p_r S^l\dP-r S^{m-1}\big(\f{P_2}{r}\p_r\dP\big)\bigg)^2dx\le C\ve^2.\eqno{(5.36)}
$$

On the other hand,
we have
$$
|S^{m-1}\mathcal{L}\dP|\le  C\ve \bigg\{\ds\sum_{0\le
l\le m}\big(r^{-\g}|S^l\p_r\dP|+r^{-\si-2}|S^lZ\dP|\big)
+\ds\sum_{0\le l\le m-1}r^{-2\g}|S^lZ^2\dP|\bigg\}.\eqno{(5.37)}
$$

Therefore, inserting (5.36)-(5.37) into (5.35) yields
\begin{align*}
&\int_{D_T}B_{2m}\dP\cdot\mathcal{M}S^m\dP dx\le C\ve^2+
C\ve\int_{D_T}
\ds\sum_{l=0}^m\bigg(r^{\mu+1-2\g+2l}(\nabla_x^l(\f{1}{r}Z\dP))^2+r^{\mu-1-\dl+2l}(\nabla_x^l\p_r\dP)^2\bigg)
dx\\
&\qquad +C\ve\bigg(\int_{D_T}r^{\mu-1-2\g}|S^{m-1}Z^2\dP|^2dx\bigg)^{\f12}.\tag{5.38}
\end{align*}

Consequently, inserting (5.9)-(5.10), (5.15)-(5.17), (5.28), (5.32)
and (5.38) into (5.6), we can complete the proof of (5.3).

For the case of $m=0$, (5.4) comes directly from Theorem 4.1, (5.10)
and (5.15)-(5.17). \qquad \qquad $\square$

Based on Lemma 5.2 and the ingredients in Lemma 5.2,
we now derive a series of estimates on the higher order derivatives of $\dP$.

{\bf Lemma 5.3. (Second order angular derivative estimates)} {\it Under the assumptions of Theorem 5.1, then
\begin{align*}
&T^{\mu}\int_{S_T}(\p_rZ\dot\Phi)^2dS
+T^{\mu-2\g}\int_{S_T}(Z^2\dot\Phi)^2dS
+\int_{D_T}\bigg(r^{\mu-1-\dl}(\p_rZ\dot\Phi)^2
+r^{\mu-1-2\g}(Z^2\dot\Phi)^2\bigg)dx\\
&\le C\ve^2+C\ve\bigg(T^{\mu}\int_{S_T}(\p_rZ\dP)^2dS+T^{\mu-2\g}\int_{S_T}(Z^2\dP)^2dS\bigg)\\
&\quad +C\ve\bigg(\ds\sum_{l=0}^1\int_{D_T}r^{\mu+1-2\g+2l}(\nabla_x^l(\f{1}{r}Z\dP))^2
+r^{\mu-1-\dl+2l}(\nabla_x^l\p_r\dP)^2dx\bigg),\tag{5.39}
\end{align*}
where $0<\dl\le\g-1$.}

{\bf Remark 5.3.} {\it Lemma 5.3, together with Lemma 5.2 for $m=1$, yields (5.2) in the case of $k=1$.}

{\bf Proof.} Noting $\mathcal {L} Z\dP=Z\mathcal {L}\dP=Z\dot f$, then it follows from Theorem 4.1 that
\begin{align*}
&T^{\mu}\int_{S_T}|\p_rZ\dot\Phi|^2dS+T^{\mu-2\g}\int_{S_T}|Z^2\dot\Phi|^2dS
+\int_{D_T}\bigg(r^{\mu-1-\dl}|\p_rZ\dot\Phi|^2+r^{\mu-1-2\g}|Z^2\dot\Phi|^2\bigg)dx\\
&\le C\int_{D_T}Z\dot f\cdot\mathcal{M}Z\dot\Phi
dx+C\ve^2.\tag{5.40}
\end{align*}

In order to estimate the term $\int_{D_T}Z\dot f\cdot\mathcal{M}Z\dot\Phi
dx$ in the right hand side of (5.40), we rewrite $Z\dot f=D_1+D_2$ with
\begin{align*}
D_1&=f_{00}\p_r^2Z\dot\Phi+\f{1}{r^2}\sum\limits_{1\le
i,j\le3}f_{ij}Z_iZ_jZ\dot\Phi
+\f{1}{r}\sum\limits_{i=1}^3f_{0i}\p_rZ_i(Z\dot\Phi),\\
D_2&=\f{1}{r^2}\sum\limits_{1\le
i,j\le3}f_{ij}[Z, Z_iZ_j]\dot\Phi
+\f{1}{r}\sum\limits_{i=1}^3f_{0i}[Z, \p_rZ_i]\dot\Phi
+Zf_{00} \p_r^2\dot\Phi\\
&\quad +Z(\f{1}{r^2}\sum\limits_{1\le i,j\le3}f_{ij}) Z_iZ_j\dot\Phi
+Z(\f{1}{r}\sum\limits_{i=1}^3f_{0i}) \p_rZ_i\dot\Phi+Zf_0,
\end{align*}
where $D_1$ contains the third order derivatives of $\dP$, and $D_2$ is composed by the lower order
(up to second order)
derivative  terms of $\dP$.

In this case, a direct computation yields
\begin{align*}
&D_1\cdot\mathcal{M}Z\dot\Phi\\
&=
\p_r\bigg(\f{1}{2}r^{\mu}a(r)f_{00}(\p_rZ\dP)^2-r^{\mu-2}a(r)\ss_{1\le i<j\le3}f_{ij}Z_iZ\dP Z_jZ\dP
-\f{1}{2}r^{\mu-2}a(r)\ss_{i=1}^3f_{ii}(Z_iZ\dP)^2\bigg)\\
&\quad +\ss_{i=1}^3Z_i\bigg(\f{1}{2}r^{\mu-1}a(r)f_{0i}(\p_rZ\dP)^2
+r^{\mu-2}a(r)\p_rZ\dP\ss_{j=1}^3f_{ij}Z_jZ\dP\bigg)\\
&\quad -\f{1}{2}\p_r(r^\mu a(r)f_{00})(\p_rZ\dP)^2-\f{1}{2}r^{\mu-1}a(r)(\p_rZ\dP)^2\ss_{i=1}^3Z_if_{0i}
+\ss_{i=1}^3\p_r(\f{1}{2}r^{\mu-2}a(r)f_{ii})(Z_iZ\dP)^2\\
&\quad +\ss_{1\le i<j\le3}\p_r(r^{\mu-2}a(r)f_{ij})Z_iZ\dP Z_jZ\dP.\tag{5.41}
\end{align*}

On the other hand, as in (5.14), it follows from the expressions of $f_{ij}, f_{0i}$ and the boundary condition
(3.4) that on $\Sigma$
\begin{align*}
&x_2\bigg(\f{1}{2}r^{\mu-1}a(r)f_{02}(\p_rZ\dP)^2+r^{\mu-2}a(r)\p_rZ\dP\ss_{j=1}^3f_{2j}Z_jZ\dP\bigg)\\
&-x_1\bigg(\f{1}{2}r^{\mu-1}a(r)f_{03}(\p_rZ\dP)^2+r^{\mu-2}a(r)\p_rZ\dP\ss_{j=1}^3f_{3j}Z_jZ\dP\bigg)=0.\tag{5.42}
\end{align*}

Thus, by integration by parts together with the expressions of $f_i$ and (5.1), we have
\begin{align*}
&|\int_{D_T}D_1\cdot\mathcal{M}Z\dP dx|\le C\ve^2+C\ve\bigg(T^{\mu}\int_{S_T}(\p_rZ\dP)^2dS
+T^{\mu-2\g}\int_{S_T}(Z^2\dP)^2dS\bigg)\\
&\qquad\qquad\qquad\qquad +C\ve\bigg(\ds\sum_{l=0}^1\int_{D_T}r^{\mu+1-2\g+2l}(\nabla_x^l(\f{1}{r}Z\dP))^2
+r^{\mu-1-\dl+2l}(\nabla_x^l\p_r\dP)^2dx\bigg).\tag{5.43}
\end{align*}

In addition, a direct computation yields
$$|D_2|\le C\ve\big(r^{-2(\g-1)}|\p_r^2\dP|+r^{-2(\g-1)}|\p_rZ\dP|+r^{-2-\si}|Z^2\dP|\big),\eqno{(5.44)}$$
which implies
$$\int_{D_T}|D_2\cdot\mathcal{M}Z\dP|dx
\le C\ve\int_{D_T}\bigg(r^{\mu+1-\dl}(\p_r^2\dP)^2+r^{\mu-1-\dl}(\p_rZ\dot\Phi)^2+r^{\mu-1-2\g}(Z^2\dP)^2\bigg)dx.\eqno{(5.45)}
$$

Substituting (5.43) and (5.45) into (5.40) yields (5.39), we then complete the proof of Lemma 5.3.
\qquad\quad $\square$

{\bf Lemma 5.4. (Third order angular derivative estimates)} {\it Under the assumptions of Theorem 5.1,
 then
\begin{align*}
&T^{\mu}\int_{S_T}(\p_rZ^2\dot\Phi)^2dS
+T^{\mu-2\g}\int_{S_T}(Z^3\dot\Phi)^2dS
+\int_{D_T}\bigg(r^{\mu-\g}(\p_rZ^2\dot\Phi)^2
+r^{\mu-1-2\g}(Z^3\dot\Phi)^2\bigg)dx\\
&\le C\ve^2+C\ve\int_{D_T}
\ss_{l=0}^2\bigg(r^{\mu+1-2\g+2l}(\nabla_x^l(\f{1}{r}Z\dP))^2+r^{\mu-1-\dl+2l}(\nabla_x^l\p_r\dP)^2\bigg)
dx\\
&\qquad +C\ve\ss_{l=0}^2\bigg(T^{\mu+2l}\int_{S_T}(\nabla_x^l\p_r\dot\Phi)^2dS
+T^{\mu-2\g+2l}\int_{S_T}(\nabla_x^l(\f{1}{r}Z\dot\Phi))^2dS\bigg),\tag{5.46}
\end{align*}
where $0<\dl\le\min\{\g-1, 2\si-2(\g-1)\}$ with $\si=\min\{1, 2(\g-1)\}$.}

{\bf Remark 5.4.} {\it Under the conditions of Lemma 5.4, as in the proof of Lemma 5.3, we have
\begin{align*}
&T^{\mu}\int_{S_T}(\p_rSZ\dot\Phi)^2dS+T^{\mu-2\g}\int_{S_T}(ZSZ\dot\Phi)^2dS\\
&\qquad +\int_{D_T}\bigg(r^{\mu-1-\dl}(\p_rSZ\dot\Phi)^2+r^{\mu-1-2\g}(ZSZ\dot\Phi)^2\bigg)dx\\
&\le C\ve^2+C\ve\int_{D_T}
\ss_{l=0}^2\bigg(r^{\mu+1-2\g+2l}(\nabla_x^l(\f{1}{r}Z\dP))^2
+r^{\mu-1-\dl+2l}(\nabla_x^l\p_r\dP)^2\bigg)
dx\\
&\qquad +C\ve\ss_{l=0}^2\bigg(T^{\mu+2l}\int_{S_T}(\nabla_x^l\p_r\dot\Phi)^2dS
+T^{\mu-2\g+2l}\int_{S_T}(\nabla_x^l(\f{1}{r}Z\dot\Phi))^2dS\bigg)\\
&\qquad +C\ve\bigg(\int_{D_T}\big(r^{\mu-1-\dl}(\p_rSZ\dP)^2+r^{\mu-1-2\g}(Z^3\dP)^2\big)dx\bigg)^{\f12}.
\end{align*}

This, together with (5.46) and (5.3) in the case of $m=2$ and Remark 2.3, yields (5.2) for $k=2$ under the
assumption that (5.2)
holds for $k\le1$.}

{\bf Proof.} At first, we establish an analogous inequality
to (5.46) for $(Z_2^2+Z_3^2)\dot\Phi$. Based on this together with
the domain decomposition technique, we can complete the proof of (5.46).

It follows from (1.9) that on $\Sigma$
$$\p_r^m\p_{\vp}\Phi=0.\eqno{(5.47)}$$

Differentiating (1.6) with respect to $\vp$ and applying (5.47) yield  on $\Sigma$
$$\p_{\vp}(\sin\vp\p_{\vp}^2\Phi)=0.\eqno{(5.48)}$$

On the other hand, by a direct computation we have
$$\p_{\vp}(Z_2^2+Z_3^2)\Phi=(\p_{\vp}^3+\cot\p_{\vp}^2-\csc^2\vp\p_{\vp})\Phi.$$

This, together with (5.47)-(5.48) and the definition of $\dot\Phi$, yields
$$
Z(Z_2^2+Z_3^2)\dot\Phi=0\q \text{on $\Sigma$}.\eqno{(5.49)}$$

Applying Theorem 4.1 to $(Z_2^2+Z_3^2)\dot\Phi$, one has for $\dl\le\g-1$
\begin{align*}
&T^{\mu}\int_{S_T}(\p_r(Z_2^2+Z_3^2)\dot\Phi)^2dS+T^{\mu-2\g}\int_{S_T}(Z(Z_2^2+Z_3^2)\dot\Phi)^2dS\\
&\quad +\int_{D_T}\bigg(r^{\mu-\g}(\p_r(Z_2^2+Z_3^2)\dot\Phi)^2+r^{\mu-1-2\g}(Z(Z_2^2+Z_3^2)\dot\Phi)^2\bigg)dx\\
&\le \int_{D_T}\mathcal{L}(Z_2^2+Z_3^2)\dot\Phi\cdot\mathcal{M}(Z_2^2+Z_3^2)\dot\Phi
dx +C\ve^2.\tag{5.50}
\end{align*}

By $\mathcal{L}(Z_2^2+Z_3^2)\dot\Phi=(Z_2^2+Z_3^2)\mathcal{L}\dot\Phi=(Z_2^2+Z_3^2)\dot f$ and the expression of $\dot f$,
we have
$$
\mathcal{L}(Z_2^2+Z_3^2)\dot\Phi=K_1+K_2+K_3+K_4\eqno{(5.51)}
$$
with
\begin{align*}
K_1&=f_{00}\p_r^2(Z_2^2+Z_3^2)\dot\Phi+\f{1}{r^2}\sum\limits_{1\le
i,j\le3}f_{ij}Z_iZ_j(Z_2^2+Z_3^2)\dot\Phi
+\f{1}{r}\sum\limits_{i=1}^3f_{0i}\p_rZ_i((Z_2^2+Z_3^2)\dot\Phi),\\
K_2&=\f{1}{r^2}\sum\limits_{1\le
i,j\le3}f_{ij}[(Z_2^2+Z_3^2),Z_iZ_j]\dot\Phi
+\f{1}{r}\sum\limits_{i=1}^3f_{0i}[(Z_2^2+Z_3^2),\p_rZ_i]\dot\Phi,\\
K_3&=\ss_{k=2}^3\sum\limits_{l=1}^2C_{l}\bigg(Z_k^{l}f_{00}Z_k^{2-l}\p_r^2\dot\Phi
+Z_k^{l}(\f{1}{r^2}\sum\limits_{1\le i,j\le3}f_{ij})Z_k^{2-l}Z_iZ_j\dot\Phi
+Z_k^{l}(\f{1}{r}\sum\limits_{i=1}^3f_{0i})Z_k^{2-l}\p_rZ_i\dot\Phi\bigg)\\
&\quad +(Z_2^2+Z_3^2)f_0^1,\\
K_4&=(Z_2^2+Z_3^2)f_0^2.
\end{align*}

Next, we start to deal with each term $\int_{D_T}K_i\cdot\mathcal{M}(Z_2^2+Z_3^2)\dot\Phi
dx$ $(1\le i\le 4)$.

\vskip 0.3 true cm
{\bf (A) The estimate on $\int_{D_T}K_1\cdot\mathcal{M}(Z_2^2+Z_3^2)\dot\Phi
dx$}
\vskip 0.3 true cm

Analogous to the treatment on $\int_{D_T}D_1\cdot\mathcal{M}(Z_2^2+Z_3^2)\dot\Phi
dx$ in (5.41) and (5.43), we can obtain
\begin{align*}
&|\int_{D_T}K_1\cdot\mathcal{M}(Z_2^2+Z_3^2)\dP dx|\\
&\le C\ve^2+C\ve\bigg(T^{\mu}\int_{S_T}(\p_r(Z_2^2+Z_3^2)\dot\Phi)^2dS+T^{\mu-2\g}\int_{S_T}(Z(Z_2^2+Z_3^2)\dot\Phi)^2dS\\
&\quad +\int_{D_T}\big(r^{\mu-1-\dl}(\p_r(Z_2^2+Z_3^2)\dot\Phi)^2+r^{\mu-1-2\g}(Z(Z_2^2+Z_3^2)\dot\Phi)^2\big)dx\bigg).
\tag{5.52}
\end{align*}

\vskip 0.3 true cm
{\bf (B) The estimate on $\int_{D_T}K_2\cdot\mathcal{M}(Z_2^2+Z_3^2)\dot\Phi
dx$}
\vskip 0.3 true cm

By the expressions of $f_{ij}$ and the assumption (5.1), it follows from a direct computation that
$$|K_2|\le C\ve
\big(r^{-2\g}|Z^3\dP|+r^{-1-\si}|\p_rZ^2\dP|\big)\eqno{(5.53)}$$
and
$$
\int_{D_T}|K_2\cdot\mathcal{M}(Z_2^2+Z_3^2)\dP|dx
\le C\ve\int_{D_T}
\ss_{l=0}^2\big(r^{\mu+1-2\g+2l}(\nabla_x^l(\f{1}{r}Z\dP))^2+r^{\mu-1-\dl+2l}(\nabla_x^l\p_r\dP)^2\big)
dx.\eqno{(5.54)}
$$

\vskip 0.3 true cm
{\bf (C) The estimate on $\int_{D_T}K_3\cdot\mathcal{M}(Z_2^2+Z_3^2)\dot\Phi
dx$}
\vskip 0.3 true cm

By the expressions of $f_{ij}$ and $f_0^1$, we know that $K_3$ only
contains such terms: $\p_r^2Z^l\dP$ ($l=0,1$), $\p_rZ^l\dP$ ($0\le
l\le 2$) and $Z^l\dP$ ($1\le l\le 3$) with suitably decayed
coefficients. More concretely, by the assumption (5.1), we have
$$|K_3|\le C\ve \big(r^{-2(\g-1)}\ss_{0\le l\le1}|\p_r^2Z^l\dP|+r^{-\g}\ss_{0\le
l\le 2}|\p_rZ^l\dP|+r^{-2\g}\ss_{1\le l\le
3}|Z^l\dP|\big),$$
which implies
$$
\int_{D_T}|K_3\cdot\mathcal{M}(Z_2^2+Z_3^2)\dP|dx
\le C\ve\int_{D_T}
\ss_{l=0}^2\bigg(r^{\mu+1-2\g+2l}(\nabla_x^l(\f{1}{r}Z\dP))^2+r^{\mu-1-\dl+2l}(\nabla_x^l\p_r\dP)^2\bigg)
dx.\eqno{(5.55)}
$$

\vskip 0.3 true cm
{\bf (D) The estimate on $\int_{D_T}K_4\cdot\mathcal{M}(Z_2^2+Z_3^2)\dot\Phi
dx$}
\vskip 0.3 true cm

By the expression of $f_0^2$, we have
$$K_4=K_4^1+K_4^2\eqno{(5.56)}$$
with
\begin{align*}
&K_4^1=\f{\hat{U}}{r^3(\hat{U}^2-c^2(\hat\rho))}\ss_{k=2}^3\ss_{i,j=1}^3\bigg(C_{ij}\big(Z_i\dP Z_k^2Z_j\dP+Z_k^2Z_i\dP Z_j\dP\big)+2Z_kC_{ij}\big(Z_i\dP Z_kZ_j\dP+Z_kZ_i\dP Z_j\dP\big)\\
&\qquad\qquad+Z_k^2C_{ij}Z_i\dP Z_j\dP\bigg)
 +\f{(\g-1)\hat{U}}{r^3(\hat{U}^2-c^2(\hat\rho))}\ss_{k=2}^3\ss_{i=1}^3Z_i\dP Z_k^2Z_i\dP,\\
&K_4^2=-\ss_{k=2}^3\bigg(\ss_{i,j=1}^3\f{2\hat{U}}{r^3(\hat{U}^2-{c}^2(\hat\rho))}C_{ij}Z_kZ_i\dP
Z_kZ_j\dP+\ss_{i=1}^3\f{2(\g-1)\hat{U}}{r^3(\hat{U}^2-c^2(\hat\rho))}(Z_kZ_i\dP)^2\bigg),\\
\end{align*}
where $K_4^1=O(\ds\f{Z\dP Z^3\dP}{r^3})+O(\ds\f{Z\dP Z^2\dP}{r^3})$, and $K_4^2=O(\ds\f{Z_iZ_j\dP Z_kZ_l\dP}{r^3})$.
Here we point out that $K_4^1$ can be easily estimated since $\ds\f{Z\dP}{r}$ admits a good decay rate in assumption (5.1).
In fact, we have
$$|K_4^1|\le
Cr^{-3}|Z\dP|\ss_{l=0}^2|Z^lZ\dP|\le Cr^{-2-\si}\ss_{l=0}^2|Z^lZ\dP|$$
and
$$
\int_{D_T}|K_4^1\cdot\mathcal{M}(Z_2^2+Z_3^2)\dP|dx
\le C\ve\int_{D_T}
\ss_{l=0}^2\big(r^{\mu+1-2\g+2l}(\nabla_x^l(\f{1}{r}Z\dP))^2+r^{\mu-1-\dl+2l}(\nabla_x^l\p_r\dP)^2\big)
dx.\eqno{(5.57)}
$$

Next we deal with $\int_{D_T}|K_4^2\cdot\mathcal{M}(Z_2^2+Z_3^2)\dP|dx$.

By H\"older inequality and Lemma 2.6-Lemma 2.7, we have for $\dl<2\si-2(\g-1)$
\begin{align*}
&|r^{\mu-3}Z_kZ_i\dP Z_kZ_j\dP\p_r(Z_2^2+Z_3^2)\dP|_{L_1(D_T)}\\
=& |r^{\f{\dl+2(\g-1)-2\si}{2}}\cdot r^{\f{2\g+2\si-9}{4}}Z_kZ_i\dP\cdot r^{\f{2\g+2\si-9}{4}}Z_kZ_j\dP
\cdot r^{\f{\mu-1-\dl}{2}}\p_r(Z_2^2+Z_3^2)\dP|_{L_1(D_T)}\\
\le &C|r^{\f{2\g+2\si-9}{4}}Z_kZ_i\dP|_{L^4(D_T)} |r^{\f{2\g+2\si-9}{4}}Z_kZ_j\dP|_{L^4(D_T)}|r^{\f{\mu-1-\dl}{2}}\p_r(Z_2^2+Z_3^2)\dP|_{L^2(D_T)}\\
\le &C\ve\ss_{l=0}^2\bigg(|r^{\f{\mu-1-\dl+2l}{2}}\nabla_x^l\p_r\dP|_{L^2(D_T)}
+|r^{\f{\mu+1-2\g+2l}{2}}\nabla_x^l(\f{1}{r}Z\dP)|_{L^2(D_T)}\bigg)
|r^{\f{\mu-1-\dl}{2}}\p_r(Z_2^2+Z_3^2)\dP|_{L^2(D_T)}\\
\le &C\ve\ss_{l=0}^2\bigg(|r^{\f{\mu-1-\dl+2l}{2}}\nabla_x^l\p_r\dP|_{L^2(D_T)}^2
+|r^{\f{\mu+1-2\g+2l}{2}}\nabla_x^l(\f{1}{r}Z\dP)|_{L^2(D_T)}^2\bigg),
\end{align*}
which derives
$$
\int_{D_T}|K_4^2\cdot\mathcal{M}(Z_2^2+Z_3^2)\dP|dx\le
C\ve\int_{D_T}
\ss_{l=0}^2\bigg(r^{\mu+1-2\g+2l}\nabla_x^l(\f{1}{r}Z\dP)^2+r^{\mu-1-\dl+2l}(\nabla_x^l\p_r\dP)^2\bigg)
dx.\eqno{(5.58)}
$$

Substituting (5.51), (5.52), (5.54)-(5.58) into (5.50) yields
\begin{align*}
&T^{\mu}\int_{S_T}(\p_r(Z_2^2+Z_3^2)\dot\Phi)^2dS
+T^{\mu-2\g}\int_{S_T}(Z(Z_2^2+Z_3^2)\dot\Phi)^2dS\\
&\qquad +\int_{D_T}\bigg(r^{\mu-\g}(\p_r(Z_2^2+Z_3^2)\dot\Phi)^2
+r^{\mu-1-2\g}(Z(Z_2^2+Z_3^2)\dot\Phi)^2\bigg)dx\\
&\le C\ve^2+C\ve\int_{D_T}
\ss_{l=0}^2\bigg(r^{\mu+1-2\g+2l}(\nabla_x^l(\f{1}{r}Z\dP))^2+r^{\mu-1-\dl+2l}(\nabla_x^l\p_r\dP)^2\bigg)
dx\\
&\qquad +C\ve\bigg(T^{\mu}\int_{S_T}(\p_r(Z_2^2+Z_3^2)\dot\Phi)^2dS
+T^{\mu-2\g}\int_{S_T}(Z(Z_2^2+Z_3^2)\dot\Phi)^2dS\bigg).\tag{5.59}
\end{align*}

By (5.59), we know that the estimates of $Z(Z_2^2+Z_3^2)\dot\Phi$ and $\p_r(Z_2^2+Z_3^2)\dot\Phi$
on $S_T$ or $D_T$ have been established but such kinds of third order angular
derivatives $Z^3\dot\Phi$ and $\p_rZ^2\dot\Phi$
are not estimated. However, away from $\vp=0$ in $\O$, we can get the estimates on the third order
derivatives $Z^3\dot\Phi$ and $\p_rZ^2\dot\Phi$. Indeed, for $\vp>\f{\vp_0}{3}$,
by Lemma 2.5 a direct computation yields
\begin{equation}
\left\{
\begin{aligned}
&\ss_{i,j,k=1}^3|Z_iZ_jZ_k\dP|\le C\big(\ss_{i=1}^3|Z_i\ss_{j=1}^3Z_j^2\dP|
+\ss_{i,j=1}^3|Z_iZ_j\dP|+\ss_{i=1}^3|Z_i\dP|\big),\\
&\ss_{i,j=1}^3|\p_rZ_iZ_j\dP|\le C\big(|\p_r\ss_{j=1}^3Z_j^2\dP|+\ss_{i=1}^3|\p_rZ_i\dP|+|\p_r\dP|\big).
\end{aligned}
\right.\tag{5.60}
\end{equation}

This, together with $Z_1^k\dot\Phi\equiv 0$ $(k\in\Bbb N)$, (5.59) and Lemma 5.3, yields
\begin{align*}
&T^{\mu}\int_{S_T\cap\{\vp\ge\f{\vp_0}{3}\}}(\p_rZ^2\dot\Phi)^2dS
+T^{\mu-2\g}\int_{S_T\cap\{\vp\ge\f{\vp_0}{3}\}}(Z^3\dot\Phi)^2dS\\
&\quad +\int_{D_T\cap\{\vp\ge\f{\vp_0}{3}\}}\bigg(r^{\mu-\g}(\p_rZ^2\dot\Phi)^2
+r^{\mu-1-2\g}(Z^3\dot\Phi)^2\bigg)dx\\
&\le C\ve^2+C\ve\int_{D_T}
\ss_{l=0}^2\bigg(r^{\mu+1-2\g+2l}(\nabla_x^l(\f{1}{r}Z\dP))^2+r^{\mu-1-\dl+2l}(\nabla_x^l\p_r\dP)^2\bigg)
dx\\
&\quad +C\ve\bigg(T^{\mu}\int_{S_T}(\p_r(Z_2^2+Z_3^2)\dot\Phi)^2dS
+T^{\mu-2\g}\int_{S_T}(Z(Z_2^2+Z_3^2)\dot\Phi)^2dS\bigg).\tag{5.61}
\end{align*}

To obtain the estimates of $Z^3\dot\Phi$ and $\p_rZ^2\dot\Phi$ on the domain $\{\vp\le\ds\f{\vp_0}{3}\}$, we will take
a domain decomposition technique. Namely, we choose
a smooth cut-off function $\chi(\vp)$ as follows
\begin{equation*}
\chi(\vp)= \left\{
\begin{aligned}
&1,\quad\text{for $0\le\vp\le\f{\vp_0}{3}$},\\
&0,\quad\text{for $\f{2\vp_0}{3}\le\vp\le\vp_0$},\\
&\text{smooth connection}, \quad \text{for $\f{\vp_0}{3}\le\vp\le\f{2\vp_0}{3}$}\\
\end{aligned}
\right.
\end{equation*}
such that $\chi(\vp)\dP$ is studied.

Indeed, by Theorem 4.1 we have
\begin{align*}
&T^{\mu}\int_{S_T}\chi(\vp)(\p_rZ^2\dot\Phi)^2dS
+T^{\mu-2\g}\int_{S_T}\chi(\vp)(Z^3\dot\Phi)^2dS\\
&\qquad +
\int_{D_T}\chi(\vp)\bigg(r^{\mu-\g}(\p_rZ^2\dot\Phi)^2
+r^{\mu-1-2\g}(Z^3\dot\Phi)^2\bigg)dx\\
&\le C\bigg(\int_{D_T}r^{\mu-1-2\g}|\chi'(\vp)|(Z^3\dot\Phi)^2dx\bigg)^{\f{1}{2}}
\bigg(\int_{D_T}r^{\mu-\g}|\chi'(\vp)|(\p_rZ^2\dot\Phi)^2dx\bigg)^{\f{1}{2}}\\
&\qquad +\int_{D_T}\mathcal{L}Z^2\dot\Phi\cdot\chi(\vp)\mathcal{M}Z^2\dot\Phi
dx+C\ve^2,\tag{5.62}
\end{align*}
here we point out that $\chi'(\vp)$ has a compact support away from $\vp=0$, which implies that
the first term in the right hand side of (5.62) can be estimated as in
(5.61).

On the other hand, by the compact support property of $\chi(\vp)$ away from $\vp=0$,
then completely similar to the treatment
on $\int_{D_T}\mathcal{L}(Z_2^2+Z_3^2)\dot\Phi\cdot\mathcal{M}(Z_2^2+Z_3^2)\dot\Phi
dx$ in (5.50), we can arrive at
\begin{align*}
&\int_{D_T}\mathcal{L}Z^2\dot\Phi\cdot\chi(\vp)\mathcal{M}Z^2\dot\Phi dx\\
\le &C\ve^2+C\ve\int_{D_T}
\ss_{l=0}^2\bigg(r^{\mu+1-2\g+2l}(\nabla_x^l(\f{1}{r}Z\dP))^2+r^{\mu-1-\dl+2l}(\nabla_x^l\p_r\dP)^2\bigg)
dx\\
&+C\ve\ss_{l=0}^2\bigg(T^{\mu+2l}\int_{S_T}(\nabla_x^l\p_r\dot\Phi)^2dS
+T^{\mu-2\g+2l}\int_{S_T}(\nabla_x^l(\f{1}{r}Z\dot\Phi))^2dS\bigg).\tag{5.63}
\end{align*}

Finally, combining (5.61) and (5.62)-(5.63) yields the proof of (5.46).\qquad\qquad $\square$

Finally, we deal with the estimates of $\na_x^4\dP$.

{\bf Lemma 5.5. (Fourth order angular derivative estimates)} {\it Under the assumptions of Theorem 5.1, then
\begin{align*}
&T^{\mu}\int_{S_T}(\p_rZ^3\dot\Phi)^2dS
+T^{\mu-2\g}\int_{S_T}(Z^4\dot\Phi)^2dS
+\int_{D_T}\bigg(r^{\mu-\g}(\p_rZ^3\dot\Phi)^2
+r^{\mu-1-2\g}(Z^4\dot\Phi)^2\bigg)dx\\
\le&C\ve^2+C\ve\int_{D_T}
\ss_{l=0}^3\bigg(r^{\mu+1-2\g+2l}(\nabla_x^l(\f{1}{r}Z\dP))^2+r^{\mu-1-\dl+2l}(\nabla_x^l\p_r\dP)^2\bigg)
dx\\
& +C\ve\ss_{l=0}^3\bigg(T^{\mu+2l}\int_{S_T}(\nabla_x^l\p_r\dot\Phi)^2dS
+T^{\mu-2\g+2l}\int_{S_T}(\nabla_x^l(\f{1}{r}Z\dot\Phi))^2dS\bigg),\tag{5.64}
\end{align*}
where $0<\dl\le\min\{\g-1, \si-(\g-1)\}$ with $\si=\min\{1, 2(\g-1)\}$.}

{\bf Remark 5.5.} {\it Under the conditions of Lemma 5.5, as in the proofs of Lemma 5.3
and Lemma 5.4 respectively, we have
\begin{align*}
&T^{\mu}\int_{S_T}(\p_rS^2Z\dot\Phi)^2dS+T^{\mu-2\g}\int_{S_T}(ZS^2Z\dot\Phi)^2dS\\
&\qquad +\int_{D_T}\bigg(r^{\mu-1-\dl}(\p_rS^2Z\dot\Phi)^2+r^{\mu-1-2\g}(ZS^2Z\dot\Phi)^2\bigg)dx\\
&\le C\ve^2+C\ve\int_{D_T}
\ss_{l=0}^3\bigg(r^{\mu+1-2\g+2l}(\nabla_x^l(\f{1}{r}Z\dP))^2
+r^{\mu-1-\dl+2l}(\nabla_x^l\p_r\dP)^2\bigg)
dx\\
&\qquad +C\ve\ss_{l=0}^3\bigg(T^{\mu+2l}\int_{S_T}(\nabla_x^l\p_r\dot\Phi)^2dS
+T^{\mu-2\g+2l}\int_{S_T}(\nabla_x^l(\f{1}{r}Z\dot\Phi))^2dS\bigg)\\
&\qquad +C\ve\bigg(\int_{D_T}\big(r^{\mu-1-\dl}(\p_rS^2Z\dP)^2+r^{\mu-1-2\g}(SZ^3\dP)^2\big)dx\bigg)^{\f12}
\end{align*}
and
\begin{align*}
&T^{\mu}\int_{S_T}(\p_rSZ^2\dot\Phi)^2dS+T^{\mu-2\g}\int_{S_T}(ZSZ^2\dot\Phi)^2dS\\
&\qquad +\int_{D_T}\bigg(r^{\mu-1-\dl}(\p_rSZ^2\dot\Phi)^2+r^{\mu-1-2\g}(ZSZ^2\dot\Phi)^2\bigg)dx\\
&\le C\ve^2+C\ve\int_{D_T}
\ss_{l=0}^3\bigg(r^{\mu+1-2\g+2l}(\nabla_x^l(\f{1}{r}Z\dP))^2
+r^{\mu-1-\dl+2l}(\nabla_x^l\p_r\dP)^2\bigg)
dx\\
&\qquad +C\ve\ss_{l=0}^3\bigg(T^{\mu+2l}\int_{S_T}(\nabla_x^l\p_r\dot\Phi)^2dS
+T^{\mu-2\g+2l}\int_{S_T}(\nabla_x^l(\f{1}{r}Z\dot\Phi))^2dS\bigg)\\
&\qquad +C\ve\bigg(\int_{D_T}\big(r^{\mu-1-\dl}(\p_rSZ^2\dP)^2+r^{\mu-1-2\g}(Z^4\dP)^2\big)dx\bigg)^{\f12}
\end{align*}

This, together with (5.64) and (5.3) in the case of $m=3$ and Remark 2.3, yields (5.2) for $k=3$.}

{\bf Proof.} As in Lemma 5.4, at first we establish an analogous inequality
to (5.53) for $Z(Z_2^2+Z_3^2)\dot\Phi$. Based on this together with
the domain decomposition technique, we can complete the proof of (5.64).

Note that
$$
\p_rZ(Z_2^2+Z_3^2)\dot\Phi=0\qquad \text{on $\vp=\vp_0$}.\eqno{(5.65)}
$$

Then applying Theorem 4.1 to $Z(Z_2^2+Z_3^2)\dot\Phi$ yields
\begin{align*}
&T^{\mu}\int_{S_T}(\p_rZ(Z_2^2+Z_3^2)\dot\Phi)^2dS
+T^{\mu-2\g}\int_{S_T}(Z^2(Z_2^2+Z_3^2)\dot\Phi)^2dS\\
&\qquad +\int_{D_T}\bigg(r^{\mu-\g}(\p_rZ(Z_2^2+Z_3^2)\dot\Phi)^2
+r^{\mu-1-2\g}(Z^2(Z_2^2+Z_3^2)\dot\Phi)^2\bigg)dx\\
&\le C\int_{D_T}\mathcal{L}Z(Z_2^2+Z_3^2)\dot\Phi\cdot\mathcal{M}Z(Z_2^2+Z_3^2)\dot\Phi
dx +C\ve^2.\tag{5.66}
\end{align*}

By $\mathcal{L}Z(Z_2^2+Z_3^2)\dot\Phi=Z(Z_2^2+Z_3^2)\mathcal{L}\dot\Phi=Z(Z_2^2+Z_3^2)\dot f$
and the expression of $\dot f$,
we have
$$
\mathcal{L}Z(Z_2^2+Z_3^2)\dP=M_1+M_2+M_3+M_4\eqno{(5.67)}
$$
with
\begin{align*}
M_1&=f_{00}\p_r^2Z(Z_2^2+Z_3^2)\dot\Phi+\f{1}{r^2}\ss_{1\le
i,j\le3}f_{ij}Z_iZ_jZ(Z_2^2+Z_3^2)\dot\Phi
+\f{1}{r}\ss_{i=1}^3f_{0i}\p_rZ_iZ(Z_2^2+Z_3^2)\dot\Phi,\\
M_2&=\f{1}{r^2}\ss_{1\le
i,j\le3}f_{ij}[Z(Z_2^2+Z_3^2),Z_iZ_j]\dot\Phi
+\f{1}{r}\ss_{i=1}^3f_{0i}[Z(Z_2^2+Z_3^2),\p_rZ_i]\dot\Phi,\\
M_3&=\ss_{k=2}^3\ss_{l=1}^2C_{l}Z\bigg(Z_k^{l}f_{00}Z_k^{2-l}\p_r^2\dot\Phi
+\ss_{1\le i,j\le3}Z_k^{l}(\f{1}{r^2}f_{ij})Z_k^{2-l}Z_iZ_j\dot\Phi
+\ss_{i=1}^3Z_k^{l}(\f{1}{r}f_{0i})Z_k^{2-l}\p_rZ_i\dot\Phi\bigg)\\
&\qquad +\bigg(Zf_{00}Z_k^2\p_r^2\dot\Phi
+\ss_{k=2}^3\ss_{1\le i,j\le3}Z(\f{1}{r^2}f_{ij})Z_k^2Z_iZ_j\dot\Phi
+\ss_{i=1}^3Z(\f{1}{r}f_{0i})Z_k^2\p_rZ_i\dot\Phi\bigg)\\
&\qquad +Z(Z_2^2+Z_3^2)f_0^1,\\
M_4&=Z(Z_2^2+Z_3^2)f_0^2.
\end{align*}
For notational simplifications, we rewrite $M_3$ and $M_4$ as follows
\begin{align*}
M_3&=\sum\limits_{l_1+l_2\le
3,l_1\ge1}C_{l_1l_2}\bigg(Z^{l_1}f_{00}Z^{l_2}\p_r^2\dot\Phi
+\sum\limits_{1\le i,j\le3}Z^{l_1}(\f{1}{r^2}f_{ij})Z^{l_2}Z_iZ_j\dot\Phi
+\sum\limits_{i=1}^3Z^{l_1}(\f{1}{r}f_{0i})Z^{l_2}\p_rZ_i\dot\Phi\bigg)\\
&\quad +Z^3f_0^1,\\
M_4&=Z^3f_0^2.
\end{align*}

Next, we treat each term $\int_{D_T}M_i\cdot\mathcal{M}Z(Z_2^2+Z_3^2)\dot\Phi
dx$ $(1\le i\le 4)$ respectively.

\vskip 0.3 true cm
{\bf (A) The estimate on $\int_{D_T}M_1\cdot\mathcal{M}Z(Z_2^2+Z_3^2)\dot\Phi
dx$}
\vskip 0.3 true cm

As in (5.52), we can obtain
\begin{align*}
&|\int_{D_T}M_1\cdot\mathcal{M}Z(Z_2^2+Z_3^2)\dP dx|\le C\ve^2+C\ve\bigg(T^{\mu}\int_{S_T}(\p_rZ(Z_2^2+Z_3^2)\dot\Phi)^2dS\\
&\qquad\quad +T^{\mu-2\g}\int_{S_T}(Z^2(Z_2^2+Z_3^2)\dot\Phi)^2dS+\int_{D_T}\big(r^{\mu-1-\dl}(\p_rZ(Z_2^2+Z_3^2)\dot\Phi)^2\\
&\qquad\quad +r^{\mu-1-2\g}(Z^2(Z_2^2+Z_3^2)\dot\Phi)^2\big)dx\bigg).
\tag{5.68}
\end{align*}

\vskip 0.3 true cm
{\bf (B) The estimate on $\int_{D_T}M_2\cdot\mathcal{M}Z(Z_2^2+Z_3^2)\dot\Phi
dx$}
\vskip 0.3 true cm
It follows from the expressions
of $f_{ij}$, (5.1) and a direct computation that
$$
|M_2|\le C\ve
\big(r^{-2\g}|Z^4\dP|+r^{-\si-1}|\p_rZ^3\dP|\big)$$
and
$$
\int_{D_T}|M_2\cdot\mathcal{M}Z(Z_2^2+Z_3^2)\dP|dx\le C\ve\int_{D_T}
\ss_{l=0}^3\bigg(r^{\mu+1-2\g+2l}(\nabla_x^l(\f{1}{r}Z\dP))^2+r^{\mu-1-\dl+2l}(\nabla_x^l\p_r\dP)^2\bigg)
dx.\eqno{(5.69)}
$$

\vskip 0.3 true cm
{\bf (C) The estimate on $\int_{D_T}M_3\cdot\mathcal{M}Z(Z_2^2+Z_3^2)\dot\Phi
dx$}
\vskip 0.3 true cm

Due to the lack of $L^{\infty}$ assumptions on the third order derivatives of $\dot\Phi$
in (5.1), we will decompose $M_3$ as follows
$$M_3=M_3^1+M_3^2,$$
where $M_3^1$ is linear with respect to the third and fourth order derivatives of $\dot\Phi$,
$M_3^2$ contains the products of two third order
derivatives of $\dot\Phi$, which is required to be specially treated and admits the following concrete
expression:
\begin{align*}
M_3^2=&-\f{C_{12}}{(\hat{U}^2-c^2(\hat\rho))}\bigg\{\bigg((\g+1)\hat{U}\p_rZ^2\dP
+(\g+1)\p_r\dP\p_rZ^2\dP+\f{\g-1}{r^2}\ss_{i=1}^3Z_i\dP
Z^2Z_i\dP\bigg)\p_r^2Z\dP\\
&-\f{\g-1}{r^2}\ss_{i=1}^3\bigg(\p_r\dP\p_rZ^2\dP-\hat{U}\p_rZ^2\dP
+\f{1}{r^2}\ss_{k=1}^3Z_k\dP Z^2Z_k\dP-\f{2}{r^2(\g-1)}Z_i\dP
Z^2Z_i\dP\bigg)ZZ_i^2\dP\\
&+\f{1}{r^4}\ss_{1\le i\neq
j\le3}Z_i\dP Z^2Z_j\dP
ZZ_iZ_j\dP+\f{1}{r^2}\ss_{i=1}^3\bigg((\hat{U}+\p_r\dP)Z^2Z_i\dP+Z_i\dP\p_rZ^2\dP\bigg)\p_rZZ_i\dP)\bigg\}
\end{align*}

By (5.1) and a direct but tedious computation, we can arrive at
$$
|M_3^1|\le C\ve \big(r^{-2(\g-1)}\ss_{0\le i \le2}|\p_r^2Z^i\dP|+r^{-\dl-1}\ss_{0\le
i\le 3}|\p_rZ^i\dP|+r^{-2-\si}\ss_{1\le i\le
4}|Z^i\dP|\big),
$$
which yields
$$
\int_{D_T}|M_3^1\cdot\mathcal{M}Z(Z_2^2+Z_3^2)\dP|dx\le C\ve\int_{D_T}
\ss_{l=0}^3\bigg(r^{\mu+1-2\g+2l}(\nabla_x^l(\f{1}{r}Z\dP))^2+r^{\mu-1-\dl+2l}(\nabla_x^l\p_r\dP)^2\bigg)
dx.\eqno{(5.70)}
$$

Next, we focus on the treatment of $\int_{D_T}|M_3^2\cdot\mathcal{M}Z(Z_2^2+Z_3^2)\dP|dx$.

By the expression of $M_3^2$, one knows that the typical terms  in $M_3^2\cdot\mathcal{M}Z(Z_2^2+Z_3^2)\dP$
are  respectively:

$r^{\mu}\p_rZ^2\dP\p_r^2Z\dP\p_r Z(Z_2^2+Z_3^2)\dP$, $r^{\mu-2}Z\dP Z^3\dP\p_r^2Z\dP\p_r Z(Z_2^2+Z_3^2)\dP$,
$r^{\mu-2}\p_rZ^2\dP Z^3\dP\p_r Z(Z_2^2+Z_3^2)\dP$,

$r^{\mu-4}Z\dP (Z^3\dP)^2\p_r Z(Z_2^2+Z_3^2)\dP$,
$r^{\mu-2}Z^3\dP\p_rZ^2\dP\p_r Z(Z_2^2+Z_3^2)\dP$ and $r^{\mu-2}Z\dP (\p_rZ^2\dP)^2\p_r Z(Z_2^2+Z_3^2)\dP$.

By the interpolation inequalities in Corollary 2.3 and Lemma 2.6, the terms mentioned above can be treated directly.
Indeed, we have

(i)
\begin{align*}
&|r^{\mu}\p_rZ^2\dP\p_r^2Z\dP\p_r Z(Z_2^2+Z_3^2)\dP|_{L^1(D_T)}\\
\le &C|r^{\dl-2(\g-1)}\cdot r^{\f{8\g-11-\dl}{4}}\p_rZ^2\dP\cdot r^{\f{8\g-7-\dl}{4}}\p_r^2Z\dP\cdot r^{\f{\mu-1-\dl}{2}}\p_rZ^3\dP|_{L_1(D_T)}\\
\le &C |r^{\f{8\g-11-\dl}{4}}\p_r^2Z\dP|_{L^4(D_T)}|r^{\f{8\g-7-\dl}{4}}\p_r^2Z\dP|_{L^4(D_T)}
|r^{\f{\mu-1-\dl}{2}}\p_rZ^3\dP|_{L^2(D_T)}.\qquad \text{(By $\dl\le\g-1$)}
\end{align*}

(ii)
\begin{align*}
&|r^{\mu-2}Z\dP Z^3\dP\p_r^2Z\dP\p_r Z(Z_2^2+Z_3^2)\dP|_{L^1(D_T)}\\
\le &|r^{\mu-1-\si}Z^3\dP\p_r^2Z\dP\p_rZ(Z_2^2+Z_3^2)\dP|_{L_1(D_T)}\qquad \text{(By assumption (5.1))}\\
\le &C|r^{\f{3\dl-(\g-1)-\si}{4}}\cdot r^{\f{4\g-11}{4}}Z^3\dP\cdot r^{\f{8\g-7-\dl}{4}}\p_r^2Z\dP\cdot r^{\f{\mu-1-\dl}{2}}\p_rZ^3\dP|_{L_1(D_T)}\\
\le &C |r^{\f{4\g-11}{4}}Z^3\dP|_{L^4(D_T)}|r^{\f{8\g-7-\dl}{4}}\p_r^2Z\dP|_{L^4(D_T)}
|r^{\f{\mu-1-\dl}{2}}\p_rZ^3\dP|_{L^2(D_T)},
\end{align*}
here we have used $0<\dl\le\min(\g-1,\si-(\g-1))$ and $\si=\min(1,2(\g-1))$, which derives $3\dl-(\g-1)-\si\le0$.

(iii)
\begin{align*}
&|r^{\mu-2}\p_rZ^2\dP Z^3\dP\p_rZ(Z_2^2+Z_3^2)\dP|_{L_1(D_T)}\\
=& |r^{\f{3\dl-4(\g-1)}{4}}\cdot r^{\f{8\g-11-\dl}{4}}\p_rZ^2\dP\cdot r^{\f{4\g-11}{4}}Z^3\dP
\cdot r^{\f{\mu-1-\dl}{2}}\p_rZ^3\dP|_{L^1(D_T)}\\
\le &C |r^{\f{8\g-11-\dl}{4}}\p_rZ^2\dP|_{L^4(D_T)} |r^{\f{4\g-11}{4}}Z^3\dP|_{L^4(D_T)}|r^{\f{\mu-1-\dl}{2}}\p_rZ^3\dP|_{L^2(D_T)}.\qquad \text{(By $\dl\le\g-1$)}
\end{align*}

(iv)
\begin{align*}
&|r^{\mu-4}Z\dP (Z^3\dP)^2\p_r Z(Z_2^2+Z_3^2)\dP|_{L^1(D_T)}\\
\le &|r^{\mu-3-\si}(Z^3\dP)^2\p_rZ(Z_2^2+Z_3^2)\dP|_{L_1(D_T)}\qquad\text{(By assumption (5.1))}\\
=&|r^{\f{\dl-2\si}{2}}\cdot (r^{\f{4\g-11}{4}}Z^3\dP)^2\cdot r^{\f{\mu-1-\dl}{2}}\p_rZ^3\dP|_{L_1(D_T)}\\
\le &C|r^{\f{4\g-11}{4}}Z^3\dP|_{L^4(D_T)}^2
|r^{\f{\mu-1-\dl}{2}}\p_rZ^3\dP|_{L^2(D_T)}.\qquad \text{(By $\dl\le \g-1\le \si$)}\\
\end{align*}

(v)
\begin{align*}
&|r^{\mu-2}Z^3\dP\p_rZ^2\dP \p_rZ(Z_2^2+Z_3^2)\dP|_{L_1(D_T)}\\
=&|r^{\f{3\dl-4(\g-1)}{4}}\cdot r^{\f{4\g-11}{4}}Z^3\dP\cdot r^{\f{8\g-11-\dl}{4}}\p_rZ^2\dP
\cdot r^{\f{\mu-1-\dl}{2}}\p_rZ^3\dP|_{L_1(D_T)}\\
\le &C|r^{\f{4\g-11}{4}}Z^3\dP|_{L^4(D_T)}|r^{\f{8\g-11-\dl}{4}}\p_rZ^2\dP|_{L^4(D_T)}
|r^{\f{\mu-1-\dl}{2}}\p_rZ^3\dP|_{L^2(D_T)}.\qquad \text{(By $\dl\le\g-1$)}\\
\end{align*}

(vi)
\begin{align*}
&|r^{\mu-2}Z\dP (\p_rZ^2\dP)^2\p_r Z(Z_2^2+Z_3^2)\dP|_{L^1(D_T)}\\
\le &|r^{\mu-1-\si}(\p_rZ^2\dP)^2\p_rZ(Z_2^2+Z_3^2)\dP|_{L_1(D_T)}\qquad\text{(By assumption (5.1))}\\
=& |r^{\dl-2(\g-1)-\si}\cdot (r^{\f{8\g-11-\dl}{4}}\p_rZ^2\dP)^2\cdot r^{\f{\mu-1-\dl}{2}}\p_rZ^3\dP|_{L_1(D_T)}\\
\le &C |r^{\f{8\g-11-\dl}{4}}\p_rZ^2\dP|_{L^4(D_T)}^2 |r^{\f{\mu-1-\dl}{2}}\p_rZ^3\dP|_{L^2(D_T)}.
\qquad \text{(By $\dl\le\g-1$)}\\
\end{align*}

Substituting those estimates above into $\int_{D_T}|M_3^2\cdot\mathcal{M}Z(Z_2^2+Z_3^2)\dP|dx$ and applying
Lemma 2.6-Lemma 2.7 yield
$$\int_{D_T}|M_3^2\cdot\mathcal{M}Z(Z_2^2+Z_3^2)\dP|dx\le
C\ve\int_{D_T}\ss_{l=0}^3\bigg(r^{\mu-1-\dl+2l}(\nabla_x^l\p_r\dot\Phi)^2
+r^{\mu+1-2\g+2l}(\nabla_x^l(\f{1}{r}Z\dot\Phi))^2\bigg)dx.\eqno{(5.71)}$$

By (5.70) and (5.71), we have
$$
\int_{D_T}|M_3\cdot\mathcal{M}Z(Z_2^2+Z_3^2)\dP|dx\le C\ve\int_{D_T}
\ss_{l=0}^3\bigg(r^{\mu+1-2\g+2l}(\nabla_x^l(\f{1}{r}Z\dP))^2+r^{\mu-1-\dl+2l}(\nabla_x^l\p_r\dP)^2\bigg)
dx.\eqno{(5.72)}
$$

\vskip 0.3 true cm
{\bf (D) The estimate on $\int_{D_T}M_4\cdot\mathcal{M}Z(Z_2^2+Z_3^2)\dot\Phi
dx$}
\vskip 0.3 true cm

Note that
\begin{align*}
M_4&=-\f{\hat{U}}{r^3}\ss_{i,j=1}^3C'_{ij}ZZ_i\dP ZZ_j\dP
-\f{\hat{U}}{r^3}\ss_{i,j=1}^3\tilde{C}_{ij}ZZ_i\dP Z^2Z_j\dP
+\f{6(\g-1)}{r^3}\hat{U}\ss_{i=1}^3ZZ_i\dP ZZ_i^2\dP\\
&\quad +\text{\{left terms\}}\\
&\equiv M_4^1+M_4^2
\end{align*}
and
$$|M_4^2|\le \f{C}{r^3}|Z\dP|\ss_{i=1}^4|Z^i\dP|.\eqno{(5.73)}$$

On the other hand, due to
\begin{align*}
&|r^{\mu-3}(Z^2\dP)^2\p_rZ(Z_2^2+Z_3^2)\dP|_{L_1(D_T)}\\
\le &|r^{\f{\dl+2(\g-1)-2\si}{2}}\cdot (r^{\f{2\g+2\si-9}{4}}Z^2\dP)^2\cdot r^{\f{\mu-1-\dl}{2}}
\p_rZ_k(Z_2^2+Z_3^2)\dP|_{L_1(D_T)}\\
\le &C|r^{\f{2\g+2\si-9}{4}}Z^2\dP|_{L^4(D_T)}^2
|r^{\f{\mu-1-\dl}{2}}\p_rZ_k(Z_2^2+Z_3^2)\dP|_{L^2(D_T)},\qquad\text{(By $\dl<2(\si-(\g-1))$)}\\
&|r^{\mu-3}Z^2\dP Z^3\dP\p_rZ(Z_2^2+Z_3^2)\dP|_{L_1(D_T)}\\
\le &|r^{\f{\dl+(\g-1)-\si}{2}}\cdot r^{\f{2\g+2\si-9}{4}}Z^2\dP\cdot r^{\f{4\g-11}{4}}Z^3\dP
\cdot r^{\f{\mu-1-\dl}{2}}\p_rZ(Z_2^2+Z_3^2)\dP|_{L_1(D_T)}\\
\le &C|r^{\f{2\g+2\si-9}{4}}Z^2\dP|_{L^4(D_T)} |r^{\f{4\g-11}{4}}Z^3\dP|_{L^4(D_T)}|r^{\f{\mu-1-\dl}{2}}\p_rZ(Z_2^2+Z_3^2)\dP|_{L^2(D_T)},
\qquad\text{(By $\dl<\si-(\g-1)$)}
\end{align*}
then together with Lemma 2.6-Lemma 2.7 we can arrive at
$$
\int_{D_T}|M_4^1\cdot\mathcal{M}Z(Z_2^2+Z_3^2)\dP|dx\le
C\ve\int_{D_T}\ss_{l=0}^3\bigg(r^{\mu-1-\dl+2l}(\nabla_x^l\p_r\dot\Phi)^2
+r^{\mu+1-2\g+2l}(\nabla_x^l(\f{1}{r}Z\dot\Phi))^2\bigg)dx.\eqno{(5.74)}$$

Collecting (5.73) and (5.74) yields
$$
\int_{D_T}|M_4\cdot\mathcal{M}Z(Z_2^2+Z_3^2)\dP|dx\le  C\ve\int_{D_T}
\ss_{l=0}^3\bigg(r^{\mu+1-2\g+2l}(\nabla_x^l(\f{1}{r}Z\dP))^2+r^{\mu-1-\dl+2l}(\nabla_x^l\p_r\dP)^2\bigg)
dx.\eqno{(5.75)}
$$

Substituting (5.67)-(5.69), (5.72) and (5.75) into (5.66) yields
\begin{align*}
&T^{\mu}\int_{S_T}(\p_rZ(Z_2^2+Z_3^2)\dot\Phi)^2dS
+T^{\mu-2\g}\int_{S_T}(Z^2(Z_2^2+Z_3^2)\dot\Phi)^2dS\\
&\qquad +\int_{D_T}\bigg(r^{\mu-1-\dl}(\p_rZ(Z_2^2+Z_3^2)\dot\Phi)^2
+r^{\mu-1-2\g}(Z^2(Z_2^2+Z_3^2)\dot\Phi)^2\bigg)dx\\
&\le C\ve^2+C\ve\int_{D_T}
\ss_{l=0}^3\bigg(r^{\mu+1-2\g+2l}(\nabla_x^l(\f{1}{r}Z\dP))^2
+r^{\mu-1-\dl+2l}(\nabla_x^l\p_r\dP)^2\bigg)
dx\\
&\qquad +C\ve\bigg(T^{\mu}\int_{S_T}(\p_rZ(Z_2^2+Z_3^2)\dot\Phi)^2dS
+T^{\mu-2\g}\int_{S_T}(Z^2(Z_2^2+Z_3^2)\dot\Phi)^2dS\bigg).\tag{5.76}
\end{align*}

By (5.76), we have obtained the estimates of $Z^2(Z_2^2+Z_3^2)\dot\Phi$ and $\p_rZ(Z_2^2+Z_3^2)\dot\Phi$
on $S_T$ or $D_T$. However, the related estimates on such third order derivatives $Z^3\dot\Phi$ and $\p_rZ^2\dot\Phi$
are not obtained. Note that away from $\vp=0$ in $\O$, as in (5.60), we can actually get the estimates of $Z^3\dot\Phi$ and $\p_rZ^2\dot\Phi$ like (5.76) since
$$
\ss_{i,j,k,l=1}^3|Z_iZ_jZ_kZ_l\dP|\le C\big(\ss_{i,k=1}^3|Z_iZ_k(Z_2^2+Z_3^2)\dP|+\ss_{i,j,k=1}^3|Z_iZ_jZ_k\dP|+\ss_{i=1}^3|Z_iZ_j\dP|+\ss_{i=1}^3|Z_i\dP|\big)$$
and
$$\ss_{i,j,k=1}^3|\p_rZ_iZ_jZ_k\dP|\le C\big(\ss_{k=1}^3|\p_rZ_k(Z_2^2+Z_3^2)\dP|+\ss_{i,j=1}^3|\p_rZ_iZ_j\dP|+\ss_{i=1}^3|\p_rZ_i\dP|+|\p_r\dP|\big)
$$
hold for $\vp>\f{\vp_0}{3}$.

Near $\vp=0$, as in (5.62), applying Theorem 4.1 yields
\begin{align*}
&T^{\mu}\int_{S_T}\chi(\vp)(\p_rZ^3\dot\Phi)^2dS
+T^{\mu-2\g}\int_{S_T}\chi(\vp)(Z^4\dot\Phi)^2dS\\
&\qquad +\int_{D_T}\bigg(r^{\mu-\g}\chi(\vp)(\p_rZ^3\dot\Phi)^2
+r^{\mu-1-2\g}\chi(\vp)(Z^4\dot\Phi)^2\bigg)dx\\
&\le C\bigg(\int_{D_T}r^{\mu-1-2\g}|\chi'(\vp)|(Z^4\dot\Phi)^2dx\bigg)^{\f{1}{2}}
\bigg(\int_{D_T}r^{\mu-\g}|\chi'(\vp)|(\p_rZ^3\dot\Phi)^2dx\bigg)^{\f{1}{2}}\\
&\qquad +\int_{D_T}\mathcal{L}Z^3\dot\Phi\cdot\chi(\vp)\mathcal{M}Z^3\dot\Phi
dx+C\ve^2,\tag{5.77}
\end{align*}
where $\chi(\vp)$ has been defined in (5.62), and
$\chi'(\vp)$ has a compact support away from $\vp=0$, which implies that
the first term in the right hand side of (5.77) can be estimated as in
(5.76).

By the compact support property of $\chi(\vp)$, similar to the treatment
on $\int_{D_T}\mathcal{L}Z(Z_2^2+Z_3^2)\dot\Phi\cdot\mathcal{M}Z(Z_2^2+Z_3^2)\dot\Phi
dx$ in (5.67), we then have
\begin{align*}
&\int_{D_T}\mathcal{L}Z^3\dot\Phi\cdot\chi(\vp)\mathcal{M}Z^3\dot\Phi dx\\
&\le C\ve^2+C\ve\int_{D_T}
\ss_{l=0}^3\bigg(r^{\mu+1-2\g+2l}(\nabla_x^l(\f{1}{r}Z\dP))^2+r^{\mu-1-\dl+2l}(\nabla_x^l\p_r\dP)^2\bigg)
dx\\
&\qquad +C\ve\ss_{l=0}^3\bigg(T^{\mu+2l}\int_{S_T}(\nabla_x^l\p_r\dot\Phi)^2dS
+T^{\mu-2\g+2l}\int_{S_T}(\nabla_x^l(\f{1}{r}Z\dot\Phi))^2dS\bigg).\tag{5.78}
\end{align*}

Substituting (5.78) into (5.77) and combining the obtained estimates for $\vp>\f{\vp_0}{3}$
in $\O$, we then complete the proof of Lemma 5.5.\qquad \qquad \qquad \qquad \qquad $\square$

Based on Lemma 5.2-Lemma 5.5 and Remark 5.1-Remark 5.5, we now start to prove Theorem 5.1.

{\bf Proof of Theorem 5.1.}

By (5.4) in Lemma 5.2, we have
$$
T^{\mu}\int_{S_T}(\p_r \dot\Phi)^2dS+T^{\mu-2\g}\int_{S_T}(Z\dot\Phi)^2dS
+\int_{D_T}\bigg(r^{\mu-1-\dl}(\p_r \dot\Phi)^2+r^{\mu-1-2\g}(Z \dot\Phi)^2\bigg)dx\le C\ve^2.\eqno{(5.79)}
$$

By $\dl\le\g-1$, Lemma 2.5 and (5.79), we can complete the proof of (5.2) in the case of $k=0$.

Similarly, by Lemma 5.2-Lemma 5.5, Remark 5.1-Remark 5.5, Lemma 2.5 and Remark 2.3, we can arrive at
\begin{align*}
&\ss_{l=1}^m\bigg(T^{\mu+2l}\int_{S_T}|\nabla_x^l\p_r\dot\Phi|^2dS
+T^{\mu-2\g+2l}\int_{S_T}|\nabla_x^l
(\f{1}{r}Z\dot\Phi)|^2dS\biggr)\\
&\quad +\ss_{l=1}^m\int_{D_T}\biggl(r^{\mu-1-\dl+2l}|\nabla_x^l\p_r\dot\Phi|^2+r^{\mu+1-2\g+2l}
|\nabla_x^l(\f{1}{r}Z\dot\Phi)|^2\biggr)dx\\
&\le C\ve^2+C\ve\biggl\{\ss_{l=0}^m\bigg(T^{\mu+2l}\int_{S_T}|\nabla_x^l\p_r\dot\Phi|^2dS
+T^{\mu-2\g+2l}\int_{S_T}|\nabla_x^l
(\f{1}{r}Z\dot\Phi)|^2dS\bigg)\\
&\quad +\ss_{l=0}^m\int_{D_T}\biggl(r^{\mu-1-\dl+2l}|\nabla_x^l\p_r\dot\Phi|^2+r^{\mu+1-2\g+2l}
|\nabla_x^l(\f{1}{r}Z\dot\Phi)|^2\biggr)dx\bigg\}\\
&\quad+C\ve\biggl(\ss_{l=1}^m\int_{D_T}\big(r^{\mu-1-\dl+2l}|\nabla_x^l\p_r\dot\Phi|^2+r^{\mu+1-2\g+2l}
|\nabla_x^l(\f{1}{r}Z\dot\Phi)|^2\big)dx\biggr)^{\f12}.\tag{5.80}
\end{align*}
Then for $1\le m\le3$, it follows from (5.79)-(5.80) that
\begin{align*}
&\ss_{l=1}^m\bigg(T^{\mu+2l}\int_{S_T}|\nabla_x^l\p_r\dot\Phi|^2dS
+T^{\mu-2\g+2l}\int_{S_T}|\nabla_x^l
(\f{1}{r}Z\dot\Phi)|^2dS\biggr)\\
&\quad +\ss_{l=1}^m\int_{D_T}\biggl(r^{\mu-1-\dl+2l}|\nabla_x^l\p_r\dot\Phi|^2+r^{\mu+1-2\g+2l}
|\nabla_x^l(\f{1}{r}Z\dot\Phi)|^2\biggr)dx\\
&\le C\ve^2+C\ve\biggl(\ss_{l=1}^m\int_{D_T}r^{\mu-1-\dl+2l}|\nabla_x^l\p_r\dot\Phi|^2+r^{\mu+1-2\g+2l}
|\nabla_x^l(\f{1}{r}Z\dot\Phi)|^2dx\biggr)^{\f12}.\tag{5.81}
\end{align*}

If
$$\ss_{l=1}^m\int_{D_T}r^{\mu-1-\dl+2l}|\nabla_x^l\p_r\dot\Phi|^2+r^{\mu+1-2\g+2l}
|\nabla_x^l(\f{1}{r}Z\dot\Phi)|^2dx\le C\ve^2,$$
then  (5.2) is derived directly;

If
$$\ss_{l=1}^m\int_{D_T}r^{\mu-1-\dl+2l}|\nabla_x^l\p_r\dot\Phi|^2+r^{\mu+1-2\g+2l}
|\nabla_x^l(\f{1}{r}Z\dot\Phi)|^2dx\ge C\ve^2,$$
then it follows from (5.81) that
\begin{align*}
&\ss_{l=1}^m\int_{D_T}\biggl(r^{\mu-1-\dl+2l}|\nabla_x^l\p_r\dot\Phi|^2+r^{\mu+1-2\g+2l}
|\nabla_x^l(\f{1}{r}Z\dot\Phi)|^2\biggr)dx\\
\le& C\ve\biggl(\ss_{l=1}^m\int_{D_T}r^{\mu-1-\dl+2l}|\nabla_x^l\p_r\dot\Phi|^2+r^{\mu+1-2\g+2l}
|\nabla_x^l(\f{1}{r}Z\dot\Phi)|^2dx\biggr)^{\f12},
\end{align*}
which means
$$\ss_{l=1}^m\int_{D_T}\biggl(r^{\mu-1-\dl+2l}|\nabla_x^l\p_r\dot\Phi|^2+r^{\mu+1-2\g+2l}
|\nabla_x^l(\f{1}{r}Z\dot\Phi)|^2\biggr)dx\le C\ve^2.$$
Substituting this into (5.81) derives (5.2) for $1\le k\le 3$ and further completes the proof of Theorem 5.1.\qquad$\square$

\vskip 0.5 true cm
{\bf $\S 6$. Proof of Theorem 1.1.}
\vskip 0.5 true cm

It follows from Sobolev's embedding theorem (see also [10, Lemma
14]) that, one has for $1\le r\le T$
\begin{equation}
\left\{
\begin{aligned}
&\sum\limits_{0\le l\le 1}|r^l\nabla_x^l(\p_r\dot\Phi)|^2\le
Cr^{-2}\int_{S_r}\sum\limits_{0\le l\le 3}
|r^l\nabla_x^l(\p_r\dot\Phi)|^2dS,\\
&\sum\limits_{0\le l\le
1}|r^l\nabla_x^l(\f{1}{r}Z\dot\vp)|^2\le
Cr^{-2}\int_{S_r}\sum\limits_{0\le l\le 3}
|r^l\nabla_x^l(\f{1}{r}Z\dot\vp)|^2dS.
\end{aligned}
\right.\tag{6.1}
\end{equation}

On the other hand, (5.2) shows that
$$\int_{S_r}\sum\limits_{0\le l\le
3}|r^l\nabla_x^l(\p_r\dot\Phi)|^2dS\le C\ve^2r^{-\mu},\q
\int_{S_r}\sum\limits_{0\le l\le
3}|r^l\nabla_x^l(\f{1}{r}Z\dot\Phi)|^2dS\le
C\ve^2r^{-\mu+2\g-2}.$$

Hence we arrive at
$$\sum\limits_{0\le l\le
1}|r^l\nabla_x^l(\p_r\dot\Phi)|^2\le C\ve^2 r^{-\mu-2},\q
\sum\limits_{0\le l\le 1}|r^l\nabla_x^l(\f{1}{r}Z\dot\Phi)|^2\le
C\ve^2r^{-\mu+2\g-4}.$$

Subsequently, one has
$$\sum\limits_{0\le l\le
1}|r^l\nabla_x^l(\p_r\dot\Phi)|\le C\ve r^{-\f{\mu+2}{2}}=C\ve
r^{-2(\g-1)},\q \sum\limits_{0\le l\le
+1}|r^l\nabla_x^l(\f{1}{r}Z\dot\Phi)|\le C\ve
r^{-\f{\mu-2\g+4}{2}}=C\ve r^{-(\g-1)}.$$

In addition, for $1<\g<2$,
$$
|Z\dP|\le|Z\dP(1,\vp)|+\int_1^r|\p_rZ\dP(t,\vp)|dt\le C\ve(1+r^{1-2(\g-1)}),
$$
which means
$$
\f{1}{r}|Z\dP|\le C\ve(r^{-1}+r^{-2(\g-1)})\le C\ve r^{-\si}.
$$

In this case, by the Bernoulli's law (1.2),
we have $c^2(\rho)=
c^2(\hat\rho)-\f{\g-1}{2}\big((\p_r\dot\Phi)^2+2\hat U\p_r\dot\Phi+\f{1}{r^2}(Z\dot\Phi)^2\big)$,
which derives $Cr^{2(1-\g)}-C\ve(r^{2(1-\g)}+r^{-2\si})<c^2(\rho)<Cr^{2(1-\g)}+C\ve(r^{2(1-\g)}+r^{-2\si})$
together with Lemma 2.1. Thus, one obtains $c^2(\rho)\sim r^{2(1-\g)}>0$ for any $r\ge 1$ and small $\ve$.
Therefore, the proof of Theorem
1.1 is completed by the local existence result in Theorem 3.1 and continuous induction method, where
the $C^{\infty}$ regularity of $\Phi$ comes from the strong continuity principle (see [19]) and the $C^{\infty}$ smoothness of initial data $(\Phi(x)|_{r=1}, \p_r\Phi(x)|_{r=1})$ and boundary $\Sigma$.\qquad\qquad $\square$

\end{document}